\documentclass[9pt]{article}
 \usepackage{float}
\usepackage{amsmath}
\usepackage{amssymb}
\usepackage{amscd}
\usepackage{hyperref}
\def\ke{k}
\newtheorem{theorem}{Theorem}[section]
\newtheorem{lemma}[theorem]{Lemma}
  \newtheorem{definition}[theorem]{Definition}
\newtheorem{corollary}[theorem]{Corollary}
\newtheorem{proposition}[theorem]{Proposition}
 \newtheorem{remark}[theorem]{Remark}
    \newtheorem{example}[theorem]{Example}

\begin{document}
 \title{Realization of a certain class of semi-groups as value semi-groups of valuations   }

\author{Mohammad Moghaddam
\\{\footnotesize School of Mathematics, Institute in Theoretical Physics and Mathematics,}
\\{\footnotesize P. O. BOX 19395-5746, TEHRAN, IRAN}
\\{\footnotesize E-mail: moghaddam@ipm.ir}
}
\date{\today}
\maketitle
\begin{abstract} Given a well-ordered  semi-group $\Gamma$ with a minimal system of generators of ordinal type at most  $\omega n$\footnote{
We use the standard notation for ordinals, where $\omega$ is the ordinal
of the positive integers, $\omega 2=\omega+\omega$ and $\omega
(t+1)=\omega t+\omega$.} and of rational rank $r$, which satisfies a positivity and increasing condition, we construct a zero-dimensional valuation centered on the ring of polynomials with $r$ variables such that the semi-group of the values of the polynomial ring is equal to $\Gamma$.
The construction uses   a generalization of Favre and Jonsson's  version of MacLane's sequence of key-polynomials \cite{FJ}.
\end{abstract}

\section{Introduction}\label{val-sec}{

Recently the interest for studying  the structure of the value semi-groups  of the valuations centered on a noetherian local-ring   has increased  (see for example \cite{CutTe}). Several examples (e.g., plane branches, irreducible quasi-ordinary hypersurface singularities) suggest that the structure of these semi-groups contains important information on the local uniformization process of the valuation. What type of semi-groups can be  realized as the semi-group of values of a  noetherian local ring dominated by a valuation ring? Little  is known in this respect. We know they are well-ordered of ordinal type $<\omega^h$, for some natural number $h$ (\cite{ZS}, Appendix 3, Proposition 2). Abhyankar's inequality holds between numerical invariants of these valuations (see below). And, such semi-groups have no accumulation point when they are considered as semi-groups of $(\mathbb{R}^n,<_{lex})$ \cite{CutTe}.

 \par In this paper we show that given a semi-group $\Gamma$ of rational rank $r$, with a given minimal system of generators which is  well-ordered of ordinal type at most $\omega n$, $n\in\mathbb{N}$, which satisfies a positivity and increasing condition (Definition \ref{positively generated} and Theorem \ref{semigroup generation}), there is a polynomial ring $R=k[X_1,\ldots,X_r]$, where $k$ is an arbitrary field,  and a valuation $\nu$, which is positive on $R$, such that the value semi-group  $\nu(R\setminus\{0\})$ is equal to $\Gamma$. %Moreover, we give a bound for $d$: $d\leq r$. It seems that this bound is optimal. More precisely, we conjecture that, given such a semi-group $\Gamma$, the least number of variables which is needed to realize $\Gamma$ as the value semi-group of a valuation centered on a polynomial ring is $n+r+1$.
 \par Our basic tool is  a generalization of Favre and Jonsson's version of MacLane's sequence of key-polynomials (\cite{FJ}, \cite{Mac1}) for polynomial rings with arbitrary number of variables. The technique of sequences of key-polynomials was first invented by
MacLane \cite{Mac1}, following ideas of Ostrowski, to produce and describe all the extensions of
a discrete rank one valuation $\nu$ of a field $K$ to the extension
field $L=K(x).$ He attached to any extension,
say $\mu,$ of the valuation $\nu,$ a  sequence of polynomials $\phi_i(x)$ of the
ring $K[x]$.  %Roughly speaking, the first key polynomial is $x$, the next one is a unitary polynomial $\phi_1(x)=x^{d_1}+a_{d_1-1}x^{d_1-1}+\cdots +a_0$ of minimal degree in the ring $K[x]$ whose value $\beta_1=\mu(\phi_1(x))$  is larger than its Gauss valuation $\mu_0(\phi_1)=\hbox{\rm inf}(\nu (a_i)+i\mu(x))$. If there is no such polynomial, the valuation $\mu$ is determined by $\mu (x)$. Otherwise, using Euclidean division to expand polynomials in the form $P(x)=\sum_jA_j(x)\phi_1(x)^j$ with degree of $A_j$ less than $d_1$ and obtain a  new valuation $\mu_1=[\mu;\mu_1(\phi_1)=\beta_1]$ determined by $\mu_1(P(x))=\hbox{\rm inf}(\mu_0 (A_j)+j\beta_1)$. Then one takes a unitary polynomial of minimal degree for which $\beta_2=\mu (\phi_2(x))>\mu_1((\phi_2(x))$ and so on.
 By induction one can produce any
extension $\mu$ to $L$ of the valuation $\nu$ using valuations
constructed by key-polynomials (augmented valuations).
 In \cite{vaq1}, Vaqui\'e generalized
MacLane's method to produce all the extensions of an arbitrary
valuation of an arbitrary field $K$ to $L.$ He showed that given
such an extension of a valuation, there may be many ways to
produce such countable well-ordered sets of key-polynomials and augmented
valuations.   Later Favre and Jonsson showed that in the case of $d=1$ one can
consider a rather simple sequence of toroidal-key-polynomials
 (SKP), to produce all the pseudo-valuations centered on
the ring $k[[X_0,X_1]]$. Using the arithmetic of the sequence of key-polynomials of the extension $\mu$ of the valuation $\nu$, in \cite{Vaq2}, Vaqui\'e defined a new invariant, called total jump (saut total). In the case where $L=K[x]$ and $x$ algebraic over $K$, he gives a formula relating total jump to the classical invariants of the valuation extensions. In \cite{Sp2}, the construction of key-polynomials is generalized for the case where $L$ is an arbitrary algebraic extension of $K$ (not necessarily of the form $K[x]$).  They give an explicit description of the construction of key-polynomials of the valuation extension $(L,\mu)$ of $(K,\nu)$.  There are several constructions in \cite{Sp2} which are analogous to the present work, for example the notion of {\sl standard monomial} and {\sl standard expansion} corresponds to the {\sl monomial of $adic$ form} and {\sl adic expansion}, respectively, in our terminology.

 In this text, we give a  generalization of
  the sequence of  toroidal-key-polynomials of \cite{FJ} to produce
   a class of valuations of the field $k((X_0,\ldots,X_d))$, where $k$ is an arbitrary field.  Our generalization cannot generate all the  valuations centered at $k[[X_0,\ldots,X_d]]$. The construction is explicit enough to describe the value semi-group $\nu(k[[X_0,\ldots,X_d]]\setminus\{0\})$. And in addition to realize certain semi-groups as value semi-groups.

Here we recall the basic definitions associated to valuations.

\begin{definition}{\rm Fix a valuation $\nu.$
\begin{itemize}
 \item The rank $\mathrm{rk}(\nu)$ of $\nu,$
   is the Krull dimension of the
valuation ring $R_{\nu}.$
\item The rational rank of $\nu,$
$\mathrm{r.rk}(\nu),$  is the dimension of
$\nu(Frac(R_{\nu})^*)\otimes_{\mathbb{Z}} \mathbb{Q}$ as a vector space
over $\mathbb{Q}.$
\item The transcendence degree of $\nu,$
$\mathrm{tr.deg}(\nu),$  is the transcendence degree of the
extension of $k$ over residue field of $\nu,$ $k\subseteq
k_{\nu}:=\frac{R_{\nu}}{\frak{m}_{\nu}}. $
\end{itemize}
 }
\end{definition}

 The principal relation between these
numerical invariants is given by Abhyankar's
inequalities:$$\mathrm{rk}(\nu)+\mathrm{tr.deg}(\nu)\leq
\mathrm{r.rk}(\nu)+\mathrm{tr.deg}(\nu)\leq \mathrm{dim}R.$$
Moreover, if $ \mathrm{r.rk}(\nu)+\mathrm{tr.deg}(\nu)=
\mathrm{dim}R,$ then value group is isomorphic (as a group) to
$\mathbb{Z}^{\mathrm{r.rk}(\nu)}.$ When $\mathrm{
rk}(\nu)+\mathrm{tr.deg}(\nu)= \mathrm{dim}R,$ the value group is
isomorphic as an ordered group to $ \mathbb{Z}^{\mathrm{
rk}(\nu)},$ endowed with the lex. order.

Let $R$ be an integral domain with field of fractions $K$ and let
$\nu$ be a valuation of $K$ such that its valuation ring $R_{\nu}$
contains $R,$ in this case we say the valuation is centered on the
ring $R.$ Let us denote by $\Phi$ the totally ordered value group
of the valuation $\nu.$ Denote by $\Phi_{+}$ the semigroup of
positive elements of $\Phi$ and set $\Gamma=\nu(R\setminus
\{0\})\subset \Phi_{+}\cup\{0\};$ it is the semigroup of
$(R,\nu);$ since $\Gamma$ generates the group $\Phi,$ it is
cofinal in the ordered set $\Phi_{+}.$

For $\phi\in\Phi,$ set
$$\mathcal{P}_{\phi}(R)=\{x\in R\mid \nu(x)\geq \phi\}$$
$$\mathcal{P}^{+}_{\phi}(R)=\{x\in R\mid \nu(x)> \phi\},$$
where we agree that $0\in \mathcal{P}_{\phi}$
 for all $\phi,$ since its value is larger than any $\phi,$ so that by the properties of valuations the $\mathcal{P}_{\phi}$
 are ideals of $R.$ Note that the intersection $\bigcap_{\phi\in
\Phi_{+}}\mathcal{P}_{\phi}=(0)$ and that if $\phi$ is in the
negative
 part $\Phi_{-}$ of $\Phi,$ then $\mathcal{P}_{\phi}(R)=\mathcal{P}^{+}_{\phi}(R)=R.$

For $\phi\notin \Gamma,$
$\mathcal{P}_{\phi}(R)=\mathcal{P}_{\phi}^{+}(R).$ For each non
zero element $x\in R,$ there is a unique $\phi\in \Gamma$ such
that $x\in \mathcal{P}_{\phi}\setminus\mathcal{P}^{+}_{\phi};$ the
image of $x$ in the quotient
$(\mathrm{gr}_{\nu}R)_{\phi}=\mathcal{P}_{\phi}/\mathcal{P}^{+}_{\phi}$
is the {\sl initial form} $\mathrm{in}_{\nu}(x)$ of $x.$

The graded algebra associated with the valuation $\nu$ was
introduced in (\cite{L-T},\cite{Te2}) for the very special case of
a plane branch (see \cite{GT}), and in \cite{Sp1} in full
generality. Later it was extensively used in \cite{Te1} as a tool
to solve the local-uniformization problem. It is
$$\mathrm{gr}_{\nu}R=\bigoplus_{\phi\in\Gamma}\mathcal{P}_{\phi}(R)/\mathcal{P}^{+}_{\phi}(R).$$

\noindent
Acknowledgments: This work has been carried out during my Ph.D.
work as a cotutelle student at Institut Mathematiques de Jussieu, in coordination with  Universite Paris Sud (Orsay) and
 University of Tehran. I am grateful to my advisor Bernard Teissier for suggesting the problem and many useful
 conversations, to Rahim Zaare-Nahandi  and Laurent Clozel  for their cooperation on this program. The author is also
 thankful to the officials of the three universities, as well as the Institute for Studies in Physics and Mathematics (IPM, Tehran)
 and the Cultural Section of the French Embassy in Tehran and Crous de Versailles who were all involved in this issue. I would like to thank Michel Waldschmidt for his help, and the referee for his thoughtful study and many valuable comments.  
}

\section{The inductive definition of SKP's}{   From now on  by $\Phi$ we  mean a totally ordered abelian group of rank $d+1$. The total ordering of $\Phi$ is denoted by $<$.
Let $\Delta_{0}=(0)\subset\cdots\subset\Delta_{d+1}=\Phi$ be its sequence
of isolated subgroups (see \cite{ZS}). We define the sequence of pre-values and  the sequence of values of positive type. associated to a sequence of values of positive type there exists   a  sequence of key-polynomials (SKP) which are  elements of the power series ring  $\ke^{(d)}=k[[X_0,\ldots,X_d]]$\footnote{for any $i\leq d$ we define $\ke^{(i)}=k[[X_0,\ldots,X_i]]$ and $\ke_{(i)}=k((X_0,\ldots,X_i))$}. First we
need a general lemma on abelian groups.

\begin{lemma}\label{halgiri}
Let $\Psi$ be an abelian group $\alpha$ an ordinal number and
$\Gamma=\{\gamma_0,\gamma_1,\ldots,\gamma_{\alpha}\}$   be a well-ordered sequence of elements
of $\Psi.$ For any ordinal $i\leq \alpha$ define the subgroups of $\Psi$,
$G_i=(\gamma_j)_{j\leq i}$\footnote{If $a_1,\ldots,a_n$ are elements of a group $G$, by $(a_1,\ldots,a_n)$ we denote the subgroup generated by these elements and by $\langle a_1,\ldots,a_n\rangle$ the semigroup generated by them.}, $G_{i^-}=(\gamma_j)_{j<i}$,  $n_i=[G_i:G_{i^-}]$, and   set $n_0=\infty.$  Then for any $i\leq \alpha$ such that $n_i\neq \infty$, we have
   a unique representation
   \begin{equation}\label{nice adic}
   n_i\gamma_i=\sum_{j<i}m_j\gamma_j,
   \end{equation} where $0\leq m_j<n_j$ when $n_j\neq \infty,$ and $m_j\in \mathbb{Z}$ when
$n_j=\infty$, and $m_j=0$ except for a finite number of $j$. More generally, every element of $G_{i^-}$ can be written uniquely in the form (\ref{nice adic}).
\end{lemma}
\noindent
{\bf Proof. } Let $i\leq \alpha$ and   $n_i\neq\infty$, by definition of $n_i$ we have $n_i\gamma_i\in G_{i^-}$. Thus, there exists a representation   $n_i\gamma_i=\sum_{j<i}p_j\gamma_j$, where $p_j\in \mathbb{Z}$, and $p_j=0$ except for a finite number of $j$.
We define, inductively,  a sequence  $A:N'\subset\mathbb{N}\to \{1,\ldots,\alpha\}$ of elements of the index set $\alpha$, as follows:

Let $j_0<i$ be the greatest ordinal number such that $n_{j_0}\neq \infty$ and $p_{j_0}\neq 0$, the ordinal $j_0$ exists- since there is only a finite number of non-zero $p_j$. Set $A(0)=j_0$. Using Euclidean division,  write $p_{j_0}=q_{j_0}n_{j_0}+r_{j_0}$, where $0\leq r_{j_0}<n_{j_0}$. Substituting this for  $p_{j_0}$, and expanding $n_{j_0}\gamma_{j_0}$ in terms of elements of $G_{j_0^-}$, we get
$n_i\gamma_i=\sum_{j<j_0}p'_j\gamma_j+r_{j_0}\gamma_{j_0}$, where $p'_j\neq 0$ except for a finite number of $j$. Now, as before, let $j_1(<j_0)$ be the first ordinal number such that $n_{j_1}\neq \infty$ and $p'_{j_1}\neq 0$. Set $A(1)=j_1$ and continue as before to obtain $n_i\gamma_i=\sum_{j<j_1}p''_j\gamma_j+r_{j_1}\gamma_{j_1}+r_{j_0}\gamma_{j_0}$, where $0\leq r_{j}<n_j$. Continue this construction.

Either this construction stops after a finite number of  steps, say $j_k$, then we have $n_i\gamma_i=\sum_{j<i}m_j\gamma_j$, such that $m_j=0$ except for a finite number of $j$, and $0\leq m_j<n_j$ when $n_j\neq\infty$. This shows the existence part of the claim in this case. Or, the construction continues for ever, in this case we get a strictly decreasing sequence $A:\mathbb{N}\to \alpha$. But this is impossible: It suffices to note that $A(\mathbb{N})$ is a subset of $\alpha$ without least element, which is impossible (as $\alpha$ is well-ordered). Thus we have proved the existence part of the claim.

For the uniqueness, if we have two such representation $n_i\gamma_i=\sum_{j<i}m_j\gamma_j=\sum_{j<i}m'_j\gamma_j$ then let $j_0$ be the greatest index such that $m_{j_0}\neq m'_{j_0}$ (as the number of nonzero $m_j$ and $m'_j$ is finite this greatest index exists). Suppose $m_{j_0}>m'_{j_0}$ then $(m_{j_0}-m'_{j_0})\gamma_{j_0}=\sum_{j<j_0}(m'_j-m_j)\gamma_j\in G_{j_0^-}$ which is a contradiction, because $0\leq m_{j_0}-m'_{j_0}<n_{j_0}$.  $\hfill\Box$\\

\begin{definition}\label{positively generated}{\rm With the notation of Lemma \ref{halgiri}, we say the sequence $\Gamma$ {\sl is of positive type in} the group $\Psi$, if for any $i$ we have all $m_j\in \mathbb{N}$.
}
\end{definition}
This positivity condition implies that for all $i$, $\gamma_i$ is in the positive cone generated by the previous $\gamma$'s. However, the converse of this is not necessarily true. This condition enable us to construct our key-polynomials as binomials in terms of previous key-polynomials (Definition \ref{SKP}).
\begin{definition} \label{sequence of positively generated values}
{\rm  A sequence  $(\beta_{i,j}\in \Phi)_{i=0..d, j=1..\tilde{\alpha}_i}$\footnote{By $i=1..d$ and $j=1..\tilde{\alpha}_i$ we mean $i=1,\ldots,d$ and $j=1,\ldots,\tilde{\alpha}_i$. },
 $\tilde{\alpha}_i$ an ordinal number and $\tilde{\alpha}_0=1$,  is
called a {\sl sequence of pre-values} if for any $i$ and $j$ we have
\begin{itemize}
\item $\beta_{i,j+1}> n_{i,j}\beta_{i,j}$, where
$n_{i,j}=\mathrm{min}\{r\in\mathbb{N}\cup\{\infty\}:\ r\beta_{i,j}\in (\beta_{i',j'})_{(i',j')<_{lex}(i,j)}\}$
 \item $n_{i,j}\neq \infty$  for $j<\tilde{\alpha}_i$.
 \item When $j$ is a limit ordinal then $\beta_{i,j}>\beta_{i,j'}$, for any $j'<j$.
\end{itemize}}
  \end{definition}

 Consider the index set $\{(i,j)\}_{i=0..d,j=1..\tilde{\alpha}_i}$, ordered by the $lex.$ ordering. As $\tilde{\alpha}_i$ are ordinals, this is a well ordering. According to Lemma \ref{halgiri}, when $n_{i,j}\neq\infty $
there exists a unique representation
\begin{equation}\label{beta representation}
n_{i,j}\beta_{i,j}=\sum_{(i',j')\in S_{i,j}\cup S^c_{i',j'}}m^{(i,j)}_{i',j'}\beta_{i',j'}.
\end{equation}
where $m^{(i,j)}_{i',j'}=0,$ except for a finite number of
$(i',j')<_{lex}(i,j)$, and $S_{i,j}=\{(i',j')\mid\ (i',j')<_{lex}(i,j),\ {m}^{(i,j)}_{i',j'}>0\},$
 $S^{c}_{i,j}=\{(i',j')\mid \ (i',j')<_{lex}(i,j),\ {m}^{(i,j)}_{i',j'}<0\}$. By Lemma
\ref{halgiri} we have $0\leq m^{(i,j)}_{i',j'}<n_{i',j'}$ if $n_{i',j'}\neq\infty$, and $m^{(i,j)}_{i',j'}\in\mathbb{Z}$ if $n_{i',j'}=\infty$. Thus, if $(i',j')\in S^{c}_{i,j}$ then $n_{i',j'}=\infty$ and, by definition of pre-values, we have $j'=\tilde{\alpha}_{i'}$.

  Let  $\Gamma=(\beta_{i,j}\in \Phi)_{i=0..d, j=1..\tilde{\alpha}_i}$, ordered by $lex $ ordering, be a sequence of pre-values. Let $\Phi_{d,\tilde{\alpha}_d}$ be the group generated by these elements. We say $\Gamma$ is a {\sl sequence of values} if it  is of positive type in $\Phi_{d,\tilde{\alpha}_d}$. This condition is equivalent to  $S^c_{i,j}=\emptyset$, for any $i$ and $j$.

\begin{definition}\label{SKP}{\rm ({\bf SKP's}) Given a sequence of values  $\Gamma=(\beta_{i,j}\in \Phi)_{i=0..d, j=1..\tilde{\alpha}_i}$, we associate  to $\Gamma$ a sequence of power series $(U_{i,j}\in \ke^{(d)})_{i=0..d, j=1..\alpha_i}$, $\alpha_i\leq \tilde{\alpha_i}$. It is  called the {\sl sequence of key-polynomials} of the sequence of  values $\Gamma$. It is defined by induction on $i$. For $i=0$, we set $\alpha_0=\tilde{\alpha}_0=1$ and $U_{0,1}=X_0$. Suppose $U_{i',j'}$ and $\alpha_{i'}$ are defined  for $i'<i$. We set $U_{i,1}=X_i$. Suppose $U_{i,j'}$ are defined for $j'<j$. Then we define $U_{i,j}$ as follows

\begin{itemize}

\item[(P$1$)] If $j$ is not a limit ordinal then
\begin{equation}\label{inductive representation key-poly}
U_{i,j}=U_{i,j-1}^{n_{i,j-1}}-\theta_{i,j-1}\prod_{(i',j')\in S_{i,j-1}}U_{i',j'}^{m^{(i,j-1)}_{i',j'}},
\end{equation}
where $\theta_{i,j}\in k^{*}$.
This can be written as
$U_{i,j}=U_{i,j-1}^{n_{i,j-1}}-\theta_{i,j-1}U^{m^{(i,j-1)}}.$
\item[(P$2$)] If $j$ is a limit ordinal then
 $$U_{i,j}=\lim_{j'\to j}U_{i,j'}\in \ke^{(i-1)}[[X_i]].$$
In Proposition  \ref{where is} we prove that this limit exists in the ring $\ke^{(i-1)}[X_i]$. If this limit is equal to zero, then we set $\alpha_i=j$, $\beta_{i,j}=\infty$,  and we stop the construction of the key-polynomials  at this step, for $i$. Otherwise,  we continue to construct $U_{i,j'}$ for $j'>j$.
\end{itemize}
If the construction of $U_{i,j}$'s continues for every $j\leq \tilde{\alpha}_i$ then we set $\alpha_i=\tilde{\alpha}_i$.

We denote an SKP by $[U_{i,j},\beta_{i,j}]_{i=0..d,j=1..\alpha_{i}}$.
}
\end{definition}

\begin{remark} \label{after def}{\rm The following remarks are in order:
\begin{itemize}
\item[(i)]
Given any SKP as above, if we consider the data
$[U_{i,j}, \beta_{i,j}]_{i=0,1,j=1..\alpha_i}$ then it is a
$\Gamma-$SKP for the ring $k[[X_0,X_1]]$ in the sense of \cite{FJ}
for the group $\Gamma=\Phi.$
\item[(ii)]
The formula of (P1) can be rewritten in the following way.
$$U_{i,j+1}=U_{i,j}^{n_{i,j}}-\theta_{i,j}U_{0}^{m^{(i,j)}_{0}}
U_{1}^{m^{(i,j)}_{1}}\cdots
U_{i-1}^{m^{(i,j)}_{i-1}}(U_{i,1}^{m^{(i,j)}_{i,1}} \cdots
U_{i,j-1}^{m^{(i,j)}_{i,j-1}}),$$ where $U_{i'}^{m^{(i,j)}_{i'}}=\prod_{j'\leq \alpha_{i'}}U_{i',j'}^{m^{(i,j)}_{i',j'}}$, for $i'=0..i-1$.

\item[(iii)] For a fixed $i$ when $\alpha_i$ is a limit ordinal:
 \begin{itemize}
 \item If there exists an infinite number of $j$ such that $n_{i,j}>1$  then we have
   \begin{itemize}
   \item $\mathrm{deg}_{X_i}(U_{i,j})\to \infty\ (j\to \alpha_i)$.
   \item We have $U_{i,\alpha_i}=\lim_{j\to \alpha_i}U_{i,j}=0$ (See Lemma \ref{SKP-to infinity}.(ii)).
   \end{itemize}
 \item Otherwise (we denote this case by writing  $U_{i,\alpha_i}\neq 0$), we have
   \begin{itemize}
    \item $n_{i,j}=1$, except for a finite number of ordinals $j$.
    \item There is some ordinal $j_0$ such that $\mathrm{deg}_{X_i}(U_{i,\alpha_i})=\mathrm{deg}_{X_i}(U_{i,j})$ and $n_{i,j}=1$, for all $j>j_0$.

    \end{itemize}
 \end{itemize}
\item[(iv)]
For any limit ordinal $j<\alpha_i$ there are only finitely many $j'<j$ such that $n_{i,j'}>1$: Suppose the contrary and let $j<\alpha$ be the an ordinal such that there is an infinitely many $j'<j$ such that $n_{i,j'}>1$. The argument of the proof of Lemma \ref{SKP-to infinity}.(ii) shows that $U_{i,j}=0$. Thus, by construction of SKP, we must have $j=\alpha_i$ which is a contradiction.
% Later we will see that the  $adic$ expansion of elements of $\ke^{(d)}$ are elements of $\ke^{(i)}$.

\item[(v)]
Given an SKP  and $d'\leq d-1$,  we have $(U_{i,j}\in \ke^{(d')})_{i=0..d', j=1..\alpha_i}$. Moreover, the data
$[U_{i,j}, \beta_{i,j}]_{i=0..d',j=1..\alpha_i}$ is an SKP  for the sequence of values $\Gamma'=(\beta_{i,j}\in \Phi)_{i=0..d', j=1..\tilde{\alpha}_i}$.
 \end{itemize} }
\end{remark}

\begin{definition}
{\rm  Let $[U_{i,j},\beta_{i,j}]_{i=0..d,j=1..\alpha_{i}}$ be  an SKP.  We define the semigroups $\Gamma_{i,j}$ and
the groups $\Phi_{i,j},$ for $i=0,\ldots, d,j=1,\ldots, \alpha_i,$ as
follows:
$$\Gamma_{i,j}=\langle\beta_{i',j'}\rangle_{(i',j')\leq_{lex}(i,j)},$$
$$\Phi_{i,j}=(\Gamma_{i,j}),$$
$$\Phi^{*}_{i,j}=\Phi_{i,j}\otimes_{\mathbb{Z}} \mathbb{Q}.$$
}
\end{definition}

\begin{definition} \label{degree-power series}{\rm
Consider a power series ring  $A=\ke^{(i)}$. The order
 of an element $M=\sum_{\mathbf{m}}c_{\mathbf{m}}X^{\mathbf{m}}$ of this ring is
$\mathrm{ord}_{A}(M)=\mathrm{ord}(M)=\mathrm{min}_{\mathbf{m},c_{\mathbf{m}\neq
0}}\{\sum_{q=0}^i \mathbf{m}_{q}\}$.
}\end{definition}

Let $[U_{i,j},\beta_{i,j}]_{i=0..d,j=1..\alpha_i}$ be an  SKP. Fix an $i\leq d$. Consider the abelian ordered group $\Phi_{i,\alpha_i}$. This group is order isomorphic to a subgroup of the ordered group $(\mathbb{R}^n,<_{lex})$, for some large enough $n$ (see \cite{Abh-Ramification-book}, Proposition 2.10).% The positivity condition of the SKP allows us to choose embedding such that $\beta_{i',j'}\in (\mathbb{R^+})^n$, for $(i',j')\leq_{lex}(i,\alpha_i)$.
\ Let us fix such an embedding and suppose $\alpha_i$ is a limit ordinal. Consider the first index $t\leq d$, such that  $\#\{(\beta_{i,j})_t\}_{1\leq j<\alpha_i}=\infty$. The index $t$ is independent of the choice of an ordered embedding of $\Phi_{i,\alpha_i}$ into $\mathbb{R}^n$; it is called the {\sl effective component for $i$ }. Notice that this $t$ exists: otherwise, we have $\#\{(\beta_{i,j})_t\}_{1\leq j<\alpha_i, t=1..n}<\infty$. On the other hand, we have $\beta_{i,1}<_{lex}\beta_{i,2}<_{lex}\cdots<_{lex}\beta_{i,\alpha_i}$. But this is impossible when all the components of $\beta_i$'s come from a finite set. Thus $t$ is well-defined. In \cite{CutTe}, it is shown that well-ordered semi-groups of ordinal type $\leq  \omega^h$, $h\in \mathbb{N}$,  have no accumulation point in
$\mathbb{R}^n$, in Euclidean topology. We show that the semi-groups of positive type  have a stronger property: The effective component of any sequence of  the elements of the semi-group tends to infinity (Lemma \ref{rank one infinity}, and Lemma \ref{generalized effective component})
\begin{proposition}\label{canonical simplicity} With the notation of the last paragraph we have:
\begin{itemize}
\item[$(i)$] There exists $ j_{(i)}$, $1\leq j_{(i)}<\alpha_i$, such that the first $(t-1)$ components of $\beta_{i,j}$ are the same (componentwise), for $j\geq j_{(i)}$, i.e., $(\beta_{i,j})_{t'}=(\beta_{i,j'})_{t'}$, for $j,j'\geq j_{(i)}$ and $t'<t$.
\item[$(ii)$] For $j>j'>j_{(i)}$ we have $(\beta_{i,j})_{t}\geq (\beta_{i,j'})_{t}$.
\item[$(iii)$] If $U_{i,\alpha_i}=0$ then:
 \begin{itemize}
 \item[$(1)$] $t=\mathrm{min}\{t'|\ 1\leq t'\leq n,\ \exists j<\alpha_i:\ (\beta_{i,j})_{t'}\neq 0\}$.
 \item[$(2)$]  $(\beta_{i,j})_{t'}=0$, for any $j<\alpha_i$ and $t'<t$.
 \item[$(3)$] $(\beta_{i,j})_{t}\to+\infty\ (j\to\alpha_i)$.
 \end{itemize}

\end{itemize}
\end{proposition}

\noindent
{\bf Proof. } The first item is a direct consequence of the definition.
For $(ii)$, by  definition of the SKP's, we have $\beta_{i,j}>_{lex}\beta_{i,j'}$. On the other hand, by $(i)$, the first $t-1$ components of $\beta_{i,j}$ and $\beta_{i,j'}$ are the same. Thus   $(\beta_{i,j})_{t}\geq (\beta_{i,j'})_{t}$.
\par For $(iii)$,  set $t_1=\mathrm{min}\{t'|\ 1\leq t'\leq n,\ \exists j<\alpha_i:\ (\beta_{i,j})_{t'}\neq 0\}$. By definition of $t_1$, we have $(\beta_{i,j})_{t'}=0$, for any $j<\alpha_i$ and $t'<t_1$. So, $t_1\leq t$. From the definition of the SKP, we deduce that  $\beta_{i,j+1}>_{lex}(\prod_{j_0\leq j'\leq j}n_{i,j'})\beta_{i,j_0}$. We choose $j_0$ such that $(\beta_{i,j_0})_{t_1}\neq 0$ (note that necessarily  $(\beta_{i,j_0})_{t_1}> 0$). As $U_{i,\alpha_i}= 0$, there is an infinite number of $j>j_0$ such that $n_{i,j}>1$ ($j\to\alpha_i$). This shows that $(\beta_{i,j})_{t_1}\to\infty\ (j\to\alpha_i)$. Thus $t=t_1$.   $\hfill\Box$

 \begin{lemma}\label{rank one infinity} Let $[U_{i,j},\beta_{i,j}]_{i=0..d,j=1..\alpha_i}$ be an  {\rm SKP}.  Fix an $i\leq d$ and let $t$ be the effective component for $i$. Suppose $\alpha_i$ is a limit ordinal then
$(\beta_{i,j})_t\to+\infty\ (j\to\alpha_{i})$.
\end{lemma}

\noindent {\bf Proof. } If $U_{i,\alpha_i}=0$, then the claim is the content of Proposition \ref{canonical simplicity}.$(iii)$.
Assume $U_{i,\alpha_i}\neq 0$. Then, by definition of $U_{i,\alpha_i}\neq 0$, there exists $j_0$ such that  $n_{i,j}=1$ for  $j>j_0$. Notice that in this case there is a finite number of $j$ (in general) such that $n_{i,j}\neq 1$ (by definition of $U_{i,\alpha_i}\neq 0$). And we have $(\beta_{i,j})_t=\sum_{(i',j')\in S_{i,j}}m^{(i,j)}_{i',j'}(\beta_{i',j'})_t$, for $j>j_0$.
 Define
$$C_i=\{ (i',j')\in S_{i,j},\ \mathrm{max}\{j_0,j_{(i)}\}\leq j<\alpha_i,\ (\beta_{i',j'})_t\neq 0\}.$$
\par If $\# C_i=\infty$ then  there exists some $i_0<i$ and an infinite number of $j'$ such that   $(i_0,j')\in C_i$, so we can speak of $j'\to\infty$.  For such $(i_0,j')$ (which are infinite in number) we have $n_{i_0,j'}> 1$,  hence $\alpha_{i_0}$ is a limit ordinal and $U_{i_0,\alpha_{i_0}}=0$.
Let $t'$ be the effective component for $i_0$. By definition of $C_i$ there is at least one $j'$ such that $(\beta_{i_0,j'})_{t}\neq 0$.  But $U_{i_0,\alpha_{i_0}}=0$, thus by Proposition \ref{canonical simplicity}.$(iii).(2)$, we have $t'= t$. As $(\beta_{i_0,j'})_{t}\to\infty\ (j'\to\infty)$, we have $(\beta_{i,j})_{t}\to\infty\ (j\to\alpha_i)$.

\par If $\#C_i<\infty$  then
$(\beta_{i,j})_t$'s are elements of the discrete lattice $L\subset \mathbb{R}$ generated by the finite set of generators
$\{(\beta_{i',j'})_t|\ (i',j')\in C_i\}$. Thus, as any bounded region of $\mathbb{R}$ contains only a finite number of elements of the lattice $L$, the sequence  $(\beta_{i,j})_t\ (j\to\alpha_i)$ cannot be contained in any bounded region of $\mathbb{R}$. On the other hand, by  Proposition \ref{canonical simplicity}.$(ii)$, this sequence is increasing, so, it goes to $+\infty$. $ \hfill\Box$

\begin{lemma}\label{SKP-to infinity}
Consider an {\rm SKP} $[U_{i,j},\beta_{i,j}]_{i=0..d,j=1..\alpha_i}$. Suppose $\alpha_i$ is a limit ordinal. Then we have the following:
\begin{itemize}
\item[$(i)$] For any
$n\in\mathbb{N}$ and $i<d$ there exists an ordinal $j^{(i)}_n$ such that
$\mathrm{ord}_{\ke^{(i-1)}[X_i]}(U^{m^{(i,j)}})>n$ for any $j>j^{(i)}_n$.
\item[$(ii)$] If  $U_{i,\alpha_i}=0$ then one can choose the above $j^{(i)}_n$ such that in addition $\mathrm{ord}_{\ke^{(i-1)}[X_i]}(U_{i,j})>n$ for any $j>j^{(i)}_n$.
\end{itemize}
\end{lemma}

\noindent {\bf Proof. } Suppose both $(i)$ and $(ii)$ are  proved for any $n'$ and $i'<i$, and also  for  $n'\leq n$ and $i$, and notice the result holds for $n=0$. We  prove them  for $n+1$ and $i$.   Suppose $t$ is the effective component  for $i$. For any vector $V\in \mathbb{R}^n$ we define $|V|$ to be its $t$th component, i.e., $|V|=(V)_t$. Let
 $$M^{*}=\mathrm{max}\{\{|\beta_{i',j'}|:\ (i',j')\in S_{i,j},\ j'\leq j^{(i')}_{n+1}\ \mathrm{when}\ i'<i, j'\leq j^{(i)}_n\ \mathrm{when}\ i'=i\}.$$
 Notice that the cardinality of this set is finite, so $M^*$ is well-defined.
\par  For $(i)$:
 \par   By Lemma \ref{rank one infinity}, we have
 $|\beta_{i,j}|\to +\infty\ (j\to
\alpha_i)$. Hence there exists $j^{(i)}_{n+1}$ such that
$|\beta_{i,j}|>(n+1)M^{*}$, for $j\geq j^{(i)}_{n+1}$. The claim is that this
number $j^{(i)}_{n+1}$ works. We can assume $j_{(i)}<j^{(i)}_{n}$ (see Proposition \ref{canonical simplicity}.$(ii)$). Suppose  $j>j^{(i)}_{n+1}$. \par  If
there exists at least one $(i,j')\in S_{i,j}$ such that $j'\geq
j^{(i)}_n$ then we are done. Indeed, if $m^{(i,j)}_{i,j'}>1$,
since $\mathrm{ord}_{\ke^{(i-1)}[X_i]}(U_{i,j'})>n$ (by induction assumption for $(ii)$, in the case $n$) then
$\mathrm{ord}_{\ke^{(i-1)}[X_i]}(U^{m^{(i,j)}})>nm^{(i,j)}_{i,j'}>n+1$.
 If $m^{(i,j)}_{i,j'}=1$, since $|\beta_{i,j'}|<|\beta_{i,j}|$ (because $n_{i,j'}>1$ and $\beta_{i,j}>_{lex}n_{i,j'}\beta_{i,j'}$, and $|.|$ preserves ordering for $j''>j_{(i)}$ ), there should be at least one
 other element
 $(i'',j'')\in S_{i,j}$. But $\mathrm{ord}_{\ke^{(i-1)}[X_i]}(U_{i'',j''})\geq 1.$ Therefore,
 we have
 $\mathrm{ord}_{\ke^{(i-1)}[X_i]}(U^{m^{(i,j)}})>\mathrm{ord}_{\ke^{(i-1)}[X_i]}(U_{i,j'})+\mathrm{ord}_{\ke^{(i-1)}[X_i]}(U_{i'',j''})>n+1.$
  \par If there exists some $(i',j')\in S_{i,j}$ such that
$i'<i$ and $j'>j^{(i')}_{n+1}$ then clearly we are done.

\par It remains the  case that for all $(i',j')\in S_{i,j}$:
\begin{itemize}
\item If $i'<i$ then $j'<j_{n+1}^{(i')}$.
\item If $i'=i$ then $j'<j_{n}^{(i)}$.
\end{itemize}
By definition of $M^{*}$ and  conditions above, we have
$|\beta_{i',j'}|<M^{*}$, for any $(i',j')\in S_{i,j}$.  Hence
$$|\beta_{i,j_{n+1}^{(i)}}|\leq|\beta_{i,j}|\leq n_{i,j}|\beta_{i,j}|=\sum_{(i',j')\in S_{i,j}}m^{(i,j)}_{i',j'}|\beta_{i',j'}|<(\sum_{(i',j')\in
S_{i,j}}m^{(i,j)}_{i',j'})M^{*}.$$
Where the first inequality holds because $|.|$ preserves ordering for $j'\geq j^{(i)}_n>j_{(i)}$ (Proposition \ref{canonical simplicity}.$(ii)$). \par But, by definition of $M^*$, we have
$|\beta_{i,j_{n+1}^{(i)}}|>(n+1)M^{*}$. Thus $n+1<\sum_{(i',j')\in
S_{i,j}}m^{(i,j)}_{i',j'}.$ Finally

$$\mathrm{ord}_{\ke^{(i-1)}[X_i]}(U^{m^{(i,j)}})\geq \sum_{(i',j')\in S_{i,j}}m^{(i,j)}_{i',j'}>n+1.$$
\par\noindent For $(ii)$:
\par   As $(i)$ holds for $n+1$ and  using induction assumption, we can find $j^{(i)}_{n+1}$ such that   $\mathrm{ord}(U_{i,j})>n,\ \mathrm{ord}(U^{m^{(i,j)}})>n+1$, for  $j>j^{(i)}_{n+1}$.  If this $j^{(i)}_{n+1}$ does not work for $(ii)$, find the first $j_0>j^{(i)}_{n+1}$ such that $n_{i,j_0}\neq 1$ (as $U_{i,\alpha_i}=0$ this $j_0$ exists) then set $j^{(i)}_{n+1}:=j_0$. It is straightforward to check that this new $j^{(i)}_{n+1}$ works also for $(ii)$. $\hfill\Box$

\begin{proposition}\label{where is}
Fix an {\rm SKP} $[U_{i,j},\beta_{i,j}]_{i=0..d,j=1..\alpha_i}$.  Then for any $(i,j)$ we have
 $U_{i,j}\in \ke^{(i-1)}[X_i]$.
\end{proposition}

\noindent {\bf Proof.} The proof is by induction on $i$ and $j$. For $i=0$
it is obvious. Suppose it is valid for indices less than $i$, we
prove it for $i$.  When $j$ is not a limit ordinal, formula
(P1) represents $U_{i,j}$ as a polynomial in terms of previous
$U$'s and the claim is obvious in this case by induction on  $j$.
 \par It remains the case when $j$ is a limit ordinal. We can assume that $j=\alpha_i$ (considering the SKP $[U_{i',j'},\beta_{i',j'}]_{i'=0..j',j'=1..\alpha'_{i'}}$, where $\alpha'_{i'}=\alpha_{i'}$ for $i'<i$ and $\alpha'_i=j$). We must
show that $\lim_{j'\to \alpha_i}U_{i,j'}\in\ke^{(i-1)}[X_i]$.   \par If  there is infinite number of $j$ such that $n_{i,j}>1$ then by Lemma \ref{SKP-to infinity}.(ii), we have $U_{i,\alpha_i}=0\in \ke^{(i-1)}[X_i]$. Thus, we can assume $n_{i,j}=1$, except for a finite number of $j$. Then by Lemma \ref{SKP-to
infinity}.(i), we have $\mathrm{ord}_{\ke^{(i-1)}[X_i]}(U^{m^{(i,j)}})\to\infty\ (j\to \alpha_i)$. By Remark \ref{after def}.(iii), we have $\mathrm{deg}_{X_i}(U^{m^{(i,j)}})$ is bounded.  Hence
$\mathrm{ord}_{\ke^{(i-1)}}(U^{m^{(i,j)}})\to\infty\ (j\to\alpha_i)$. Using
this fact and the equality
$U_{i,j+1}-U_{i,j}=-\theta_{i,j}U^{m^{(i,j)}}$, for $j\geq j_0$
(where $n_{i,j}=1$, for $j\geq j_0$),
 we have $\lim_{j\to\alpha_i}U_{i,j}=U_{i,j_0}^{n_{i,j_0}}-\sum_{j,j_0\leq
j<\alpha_i}\theta_{i,j}U^{m^{(i,j)}}\in\ke^{(i-1)}[X_i]$.
$\hfill\Box$

\begin{remark}\label{limit SKP}{\rm The proof of the proposition shows that for any two ordinals $j'<j''$ such that
$n_{i,j}=1$, for $j'< j< j''$, we have $U_{i,j''}=\lim_{j\to
j''}U_{i,j}=U_{i,j'}^{n_{i,j'}}-\sum_{j,j'\leq
j<j''}\theta_{i,j}U^{m^{(i,j)}}$. }
\end{remark}

%\noindent
%\begin{remark}{\rm One can see that given any SKP with $\beta$'s in a rank one ordered group, the sequence $\beta_{1,j}$
%is always bounded ($U_{1,\alpha_1}\neq 1$).  Thus it is not possible to fulfill   the assumption of the lemma in the case $d=1$. %However,  this is not true for arbitrary $d\geq 2$. On the other hand, the condition of the lemma, is a weakest possible (See %Example \ref{bounded SKP}).}
%\end{remark}

%\noindent
%\begin{example}{\rm Consider the ring $k[X_0,X_1,X_2]$ and the
%group $\Phi=\mathbb{Z}^3$ with reverse lexicographical order.
%Consider the valuation $\nu$ centered on this ring defined by the
%SKP $$[(U_{0,1},U_{1,1},U_{2,1},\ldots,U_{2,j},\ldots,U_{2,\omega}),((1,0,0),(1,1,0),(0,2,0),\ldots,(0,j+1,0),\ldots,(0,0,3))],$$
%i.e., we have $\beta_{0,1}=(1,0,0),\ \beta_{1,1}=(1,1,0),\
%\beta_{2,1}=(0,2,0),\ldots,
%\beta_{2,i}=(0,j+1,0),\ldots,\beta_{2,\omega}=(0,0,3).$ We see
%that
%$$(j+1)\beta_{0,1}+\beta_{2,j}=(j+1)\beta_{1,1}.$$ This relation
%shows that the restriction on the choice of $\beta$'s is
%satisfied. By choosing the $\theta$'s all be $1,$ the
%key-polynomials are defined as follows
%$$U_{2,j+1}=U_{2,j}-U_{0,1}^{-j-1}U_{1,1}^{j+1},\  \mathrm{for} \ j\geq 1.$$
%This example shows that it is possible to have
%$S^c_{i,j}\neq\emptyset $ for a fixed $i$ and infinite number of
%$j$'s and at the same time $n_{i,j}=1.$ }
% \end{example}

\noindent
\begin{example}{\rm

Consider the ring $k[X_0,X_1,X_2]$ and the group
$\Phi=\mathbb{Z}^3$ with reverse lexicographical order. Consider
the valuation $\nu$ centered on this ring defined by the SKP
$(U_{0,1},U_{1,1},(U_{2,j})_{j=1}^{ \omega^2})$ and
$\beta_{0,1}=(1,0,0), \beta_{1,1}=(0,1,0), \beta_{2,\omega
n+j}=(j,n+2,0)$ for $n\in\mathbb{N}, 0< j<\omega$ and
$\beta_{2,\omega^2}=(0,0,1).$ Here we have the relations
$$U_{2,\omega n+j+1}=U_{2,\omega n+j}-U_{0,1}^{j}U_{1,1}^{n+2}.$$
In this example we have $n_{2,j}=1$ for any $1<j<\omega^2.$ We see
that we cannot continue to define $U_{2,\omega^2+1}:$ the reason
is that $(\beta_{2,\omega n})_2=n+2\to \infty \ (n\to \infty)$
 and therefore necessarily
$\beta_{2,\omega^2}\notin \mathbb{Z}^2\oplus\{0\}$. Thus, as $\beta_{0,1},\beta_{1,1}\in
\mathbb{Z}^2\oplus \{0\}$ there does not exist any relation
between $\beta_{2,\omega^2},\beta_{0,1},\beta_{1,1}$ and we are
forced to stop at this step.
}
\end{example}

\noindent
\begin{example}
{\rm Consider the ring $k[X_0,X_1,X_2]$ and the group
$\Phi=\mathbb{Q}$ with the usual order $\leq.$ Consider the
valuation $\nu$ centered on this ring  by the SKP\\ $(U_{0,1},
(U_{1,j})_{j=1}^{ \omega},(U_{2,j})_{j=1}^{ \omega},\beta_{i,j})$
which is defined as follows:  Let $\{p_i\}_{i=1}^{\infty} $ be an increasing sequence of prime numbers. Define $\beta_{0,1}=1,$
$\beta_{1,1}=\frac{1}{p_1},$ $\beta_{1,j}=m_j+\frac{1}{p_j},$ for
$j\geq 2$ where $m_2=1$ and $m_{j+1}=p_jm_j+1,$ and
$\beta_{2,j}=\beta_{1,j},$ for $j\geq 1.$ Then after setting
$\theta_{i,j}=1,$  we have
$U_{1,j+1}=U_{1,j}^{p_j}-U_{0,1}^{m_{j+1}}$ and
$U_{2,j+1}=U_{2,j}-U_{1,j}.$

}
\end{example}

\section{\textit{adic} expansions}{

Suppose given an {\rm SKP} $[U_{i,j},\beta_{i,j}]_{i=0..d,j=1..\alpha_i}$. In this section we show that any element $f$ of
the power series ring $\ke^{(d)}$ has a unique expansion
in terms of key-polynomials. We give an algorithm for
computing this expansion. The algorithm is based on the notion of acceptable vectors $\alpha'\leq \alpha$ associated to the SKP. Any acceptable vector determines an SKP  $[U_{i,j},\beta_{i,j}]_{i=0..d,j=1..\alpha'_i}$. We define the notion of $(U)_{\alpha'}-adic$ expansion and show how one can get $(U)_{\alpha''}-adic$ expansions for $\alpha''\geq \alpha'$, using $(U)_{\alpha'}-adic$ expansion. In the next section, we use the $adic$ expansion of the elements to define a valuation, associated to a given SKP.

\begin{lemma}\label{SKP-poly}
Fix an {\rm SKP} $[U_{i,j},\beta_{i,j}]_{i=0..d,j=1..\alpha_i}$. When $U_{i,j}\neq 0$ it is of the form
$$U_{i,j}=X_i^{d_{i,j}}+a_{i,j,d_{i,j}-1}X_i^{d_{i,j}-1}+\cdots+a_{i,j,0}$$
where  $a_{i,j,j'}\in\ke^{(i-1)}$, such that the constant term of
$a_{i,j,j'}$ is zero. Moreover, when $j$ is not a limit
ordinal, we have $d_{i,j}=n_{i,j-1}d_{i,j-1}$ for $1\leq j< \alpha_i.$ If
$j$ is a limit ordinal then there exists an ordinal  $j_0<j$,
which is not a limit ordinal and for any $j'$ such that $j_0\leq
j'\leq j$, we have $d_{i,j'}=d_{i,j_0}=n_{i,j_0-1}d_{i,j_0-1}.$
\end{lemma}

\noindent {\bf Proof. } The proofs are all by induction. We prove
the last part. By definition of SKP's, it is clear that for any
$j'=1,\ldots, j-1,$ we have $m_{i,j'}^{(i,j)}\in S_{i,j},$ so we
have $0\leq m_{i,j'}^{(i,j)}<n_{i,j'}.$ By induction we have
$n_{i,j'}=d_{i,j'+1}/d_{i,j'}.$  Hence $m_{i,j'}^{(i,j)}+1\leq
d_{i,j'+1}/d_{i,j'}.$ So we have
$$\sum_{j'=1}^{j-1}m_{i,j'}^{(i,j)}d_{i,j'}\leq \sum_{j'=1}^{j-1}(\frac{d_{i,j'+1}}{d_{i,j'}}-1)d_{i,j'}=d_{i,j}-1<n_{i,j}d_{i,j}.$$
Hence $\mathrm{deg}_{X_i}(U_{i,j+1})=n_{i,j}d_{i,j}.$ For the last claim
 we note that when $j$ is a limit ordinal there exists a $j_0$ such that for any $j'$, $j_0\leq
j'\leq j$,  we have $n_{i,j'}=1$. $\hfill \Box$

\begin{lemma} \label{deg}
For any {\rm SKP} $[U_{i,j},\beta_{i,j}]_{i=0..d,j=1..\alpha_i}$, if $U_{i,j}\neq 0$ we have
$$\mathrm{deg}_{X_i}(U_{i,j})>\mathrm{deg}_{X_i}(\prod_{j'<j}U_{i,j'}^{p_{i,j'}}),$$ when $0\leq p_{i,j'}<n_{i,j'}.$
 In other words $\sum_{j'<j}p_{i,j'}d_{i,j'}<d_{i,j}.$ Notice that $p_{i,j'}=0$, except for a finite number of $j'$.\\
\end{lemma}

\begin{definition}
{\rm Fix an SKP $[U_{i,j},\beta_{i,j}]_{i=0..d,j=1..\alpha_i}.$ We say that a vector $(\alpha '_0,\ldots,\alpha '_d)$ such that
$\alpha '_i\leq \alpha_i$ is an acceptable vector
  if for any $i=0,\ldots, d$ and any $j=0,\ldots, {\alpha '}_{i}$ and for any $(i',j')\in S_{i,j}$ we have
$(i',j')\leq_{lex}(i',\alpha_{i'})$ for $i'<i$, and $(i',j')<_{lex}(i,j)$ when $i'=i$.
This means that in the equation (P1) defining
$U_{i,j}$ in terms of $U$'s with smaller indices, one needs only
indices from $\alpha '$, not necessarily all of $\alpha.$ Notice that an
acceptable vector $\alpha '$ determines an SKP, i.e.,
$[U_{i,j},\beta_{i,j}]_{i=0..d,j=1..\alpha '_i}$ is an
SKP.}
\end{definition}

Given an SKP $[U_{i,j},\beta_{i,j}]_{i=0..d,j=1..\alpha_i}$ the vector $\alpha$ is an acceptable vector. Moreover, the vector $(1,\ldots,1)\in\mathbb{N}^d$ is an acceptable vector for an arbitrary SKP.

\begin{definition}{\rm
Given any SKP and any acceptable $\alpha '$, one can consider the new
SKP defined by this acceptable vector and construct the power series ring $k_{((\alpha',i))}=k[[(U_{i',j'})_{i'\leq i, j'<\alpha'_i, n_{i',j'}\neq 1},(U_{i',\alpha'_{i'}})_{i'\leq i}]]\subseteq k^{(d)}$. We have $\ke^{(i)}=k_{((\alpha,i))}$.
}\end{definition}

Given an SKP $[U_{i,j},\beta_{i,j}]_{i=0..d,j=1..\alpha_i}$ and an acceptable vector $\alpha '=(\alpha '_0,\ldots,\alpha '_d),$ we want to
expand an arbitrary element  $f\in \ke^{(d)}$ in terms
of $U$'s   as an element of the power series ring $k_{((\alpha ',d))}.$

\begin{definition}\label{adic-expansion def}
  {\bf ($adic$ expansions)} \rm{Fix  an {\rm SKP} $[U_{i,j},\beta_{i,j}]_{i=0..d,j=1..\alpha_i} $.  Let $\alpha '$ be an acceptable vector for this {\rm SKP}.
For an element $f\in \ke^{(d)}$ consider the
  expansion $f=\sum_{I(J)}c_{I(J)}U^{I(J)}\in k_{((\alpha ',d))},$  where $I(J)\in\mathbb{N}^{1}\times\cdots
  \times\mathbb{N}^{\alpha '_i}\times\cdots\times\mathbb{N}^{\alpha '_d} $,   and $c_{I(J)}\in k$. This expansion is called the
   $(U)_{\alpha '}-adic$
  expansion of $f,$ when for every monomial $U^{I(J)}$  we have $ 0\leq I(J)_{i,j} < n_{i,j}$, for any  $0\leq j < \alpha'_i$
and $i=0,\ldots,d$. Notice that $I(J)_{i,j}=0,$ except for a finite number of $j$.

}
\end{definition}

\begin{definition}{\rm
 Fix  an {\rm SKP} $[U_{i,j},\beta_{i,j}]_{i=0..d,j=1..\alpha_i}$ and let $\alpha ' $ be an acceptable vector. For any monomial
$M(U)=U^{\mathbf{a}}\in k_{((\alpha ',d))}$, we define
$$\mathrm{Vdeg}(M)=(\mathrm{deg}_{X_0}(U_{0}^{\mathbf{a}_0}),\mathrm{deg}_{X_1}(U_{1}^{\mathbf{a}_1}),
\ldots,\mathrm{deg}_{X_d}(U_{d}^{\mathbf{a}_d}))\in\mathbb{N}^d.$$

%Notice that in general $U_i^{\mathbf{a}_i}\in \ke^{(i-1)}(X_i)$ and if $\mathbf{a}_{i,\alpha '_i}<0$ then $\mathrm{deg}_{X_i}(U_i^{\mathbf{a}_i})=\sum_{j=1}^{\alpha_i-1}\mathrm{deg}_{X_i}(U_{i,j}^{\mathbf{a}_{i,j}})- \mathrm{deg}_{X_i}(U_{i,\alpha '_i}^{\mathbf{a}_{i,\alpha '_i}})<0$
}
\end{definition}

\begin{definition} {\rm  Fix  an {\rm
SKP} $[U_{i,j},\beta_{i,j}]_{i=0..d,j=1..\alpha_i}$ and let $\alpha ' $ be an acceptable vector. Let
$M(U)=cU^{\mathbf{a}}$ be a monomial of the ring $k_{((\alpha ',d))}$ we
say that it is a monomial of $adic$ form if it satisfies the
conditions of monomials  of Definition \ref{adic-expansion def}. }
\end{definition}

\begin{lemma} \label{deg-rev}Fix  an {\rm SKP} $[U_{i,j},\beta_{i,j}]_{i=0..d,j=1..\alpha_i}$ and let $\alpha '$ be an acceptable vector.
 Let $M(U)=cU^{\mathbf{a}}\in k_{((\alpha ',d))}$ be a monomial of adic form with respect to this
 {\rm SKP.} Then  $\mathrm{Vdeg}(M)$ determines the vector $\mathbf{a}.$
\end{lemma}

\noindent {\bf Proof. } This is a simple consequence of Lemma
\ref{deg}. If we set $n=\mathrm{deg}_{X_i}(U_{i}^{\mathbf{a}_i})$
then we have $\mathbf{a}_{i,\alpha '_i}=[\frac{n}{d_{i,\alpha '_i}}].$ Suppose
by induction we obtained
$\mathbf{a}_{i,\alpha '_i},\ldots,\mathbf{a}_{i,j+1}$ then we have:
$\mathbf{a}_{i,j}=[\frac{n-\sum_{j'=j+1}^{\alpha '_i}\mathbf{a}_{i,j'}.d_{i,j'}}{d_{i,j}}]$. Note that if $\mathbf{a}_{i,j}\neq 0$ then for any
$j'<j$ such that $d_{i,j}=d_{i,j'}$ we have $\mathbf{a}_{i,j'}=0.$
This shows that in the case of $\alpha '_i$ be of infinite ordinal type
also the number of entries of $\mathbf{a}$ computed inductively
above, which are nonzero is finite. $\hfill \Box$

\begin{corollary}\label{well ordering} Fix  an {\rm SKP} $[U_{i,j},\beta_{i,j}]_{i=0..d,j=1..\alpha_i} $.  Let $\alpha '$ be an acceptable vector.
 For any two different monomials $M,M'$  of the power series ring $k_{((\alpha ',d))}$, we say $M<M'$ if
$$\mathrm{Vdeg}(M)<_{lex}\mathrm{Vdeg}(M').$$ This  is a
well ordering on the set of monomials of $k_{((\alpha,d))}$ of $adic$
form.
 \end{corollary}

The following proposition shows that the $adic$ expansions are
well defined elements of the ring $k_{((\alpha,d))}$ and they are unique and it gives an
algorithm to compute them.

\begin{proposition}\label{alg} {\rm \bf ( Algorithm for getting $adic$ expansions)}
 Fix an {\rm SKP} $[U_{i,j},\beta_{i,j}]_{i=0..d,j=1..\alpha_i}.$ Let $\alpha '$ and $\alpha ''$ be two acceptable vectors for this
 {\rm SKP} such
  that $\alpha '<\alpha '',$ with respect to the partial product order of $\mathbb{Z}^{d+1}.$
Let $f\in \ke^{(d)}$ and suppose we know its
${(U)}_{\alpha '}-adic$ expansion. In order to obtain its
${(U)}_{\alpha ''}-adic$ expansion we do the following:

Starting from  $(U)_{\alpha '}-adic$ expansion of $f,$ for any monomial
$M(U)$ in the expansion, and for any $i=0,\ldots,d$ and $ j<
\alpha ''_i,$ do the following replacements, and iterate this process on
the resulting expansion as far as possible.
\begin{itemize}
\item   If  $n_{i,j+1}>1$,   replace any occurrence of $U_{i,j}^{n_{i,j}}$
 in $M(U)$ by  $U_{i,j+1}+\theta_{i,j}U^{m^{(i,j)}}$  (cf. $(P1)$ of Definition \ref{SKP}).

\item If $n_{i,j+1}=1$ then let $j+1<j_0\leq \alpha ''_i$ be the first
ordinal such that $n_{i,j_0}>1$ or $j_0=\alpha ''_i$ and replace any
occurrence of $U_{i,j}^{n_{i,j}}$ in $M(U)$ by
$$U_{i,j}^{n_{i,j}}=U_{i,j_0}+\sum_{j\leq j'<j_0}\theta_{i,j'}U^{m^{(i,j')}},$$
(cf. Remark \ref{limit SKP}).
\end{itemize}
The resulting expansion is equal to the $(U)_{\alpha ''}-adic$ expansion
of the element $f.$ Moreover, this expansion is unique.
\end{proposition}

\noindent {\bf Proof. }  For any element of $k_{((\alpha,d))}$ we define
$\mathcal{M}_n$ to be those monomials with $\mathrm{ord}= n$. By
Lemma $\ref{SKP-to infinity}$, we know that $\# \mathcal{M}_n$ is
finite. We do the replacements of the algorithm (staring from
$\alpha '-adic$ expansion of $f$) in the $n-$th step only on the
monomials of $\bigcup_{n'\leq n}\mathcal{M}_{n'}$ of the current
expansion. Using Lemma $6.6$ of \cite{Mog2}, this process
terminates after finitely many steps. In this step all the
monomials of $\bigcup_{n'\leq n}\mathcal{M}_{n'}$ of the current
expansion are of $\alpha ''-adic$ form. Moreover, there exists a number
$m(n)< n$, where $m(n)\to\infty\ (n\to\infty)$, such that
in  the process of replacements on the monomials of
$\bigcup_{n'\leq n}\mathcal{M}_{n'}$ the monomials of
$\bigcup_{m'\leq m(n)}\mathcal{M}_{m'}$ does not change (Lemma
\ref{SKP-to infinity}). Doing this process as $n\to\infty$ we get
an expansion, which satisfies all the properties of $\alpha ''-adic$
expansion. Thus we obtain a $(U)_{\alpha ''}-adic$ expansion of $f$.

Now, we prove that this expansion is unique. Suppose an element
$f\in \ke^{(d)}$ has two different $adic$ expansions
 $f=\sum_{I(J)}c_{I(J)}U^{I(J)}=\sum_{I''(J'')}c''_{I''(J'')}U^{I''(J'')}.$
Assume by induction on $d$ the claim is valid for the power series ring $R\otimes_k \ke^{(d-1)}$, where
$R=k[[(U_{d,j})_{j<\alpha_d,n_{d,j}\neq 1},U_{d,\alpha_d}]]$ is
considered as the coefficient ring.
 Consider  $f$ as an
element of the ring $R\otimes_k \ke^{(d-1)}$. The two $adic$
expansions of $f$ give two $adic$ expansion of $f\in R\otimes_k
\ke^{(d-1)}$ as follows. Setting $U=(U_{(d-1)},U_{d})$ and
$I(J)=(I(J)_{(d-1)},I(J)_d)$ we have
\begin{center}
\begin{math}
\begin{array}{ll}f&=\sum_{I(J)_{(d-1)}}(\sum_{I'(J')_d,I(J)_{(d-1)}=I'(J')_{(d-1)}}c_{I'(J')}U_{d}^{I'(J')_d})U_{(d-1)}^{I(J)_{(d-1)}}\\
&=\sum_{I''(J'')_{(d-1)}}(\sum_{I'(J')_d,I''(J'')_{(d-1)}=I'(J')_{(d-1)}}c''_{I'(J')}U_{d}^{I'(J')_d})U_{(d-1)}^{I''(J'')_{(d-1)}}.
\end{array}
\end{math}
\end{center}
By induction hypothesis, these two $adic$ expansions are the same.
Suppose $M$ be the least monomial of this expansion, with respect
to the ordering of Corollary \ref{well ordering}, which refers to
the indices $I_0(J_0)$ and $I''_0(J''_0)$ (respectively). Then
equating the coefficient of $M$ in two $adic$ expansions we have

$$g=\sum_{I'(J')_d,I_0(J_0)_{(d-1)}=I'(J')_{(d-1)}}c_{I'(J')}U_{d}^{I'(J')_d}
=\sum_{I'(J')_d,I''_0(J''_0)_{(d-1)}=I'(J')_{(d-1)}}c''_{I'(J')}U_{d}^{I'(J')_d}.$$
Write $g\mid_{X_0=0,\ldots,X_{d-1}=0}=\sum_{\alpha\in
\mathbb{Z}}c_{\alpha}X_d^{\alpha}$. Let $\alpha_0$ be
the first $\alpha$ such that $c_{\alpha}\neq 0$. Then by Lemma
\ref{SKP-poly} and \ref{deg-rev}, there is a unique monomial in
either of the expansions of $g$ ($M$ and $M'$ respectively) such
that $\mathrm{Vdeg}(M)=\mathrm{Vdeg}(M')=\alpha$ (Here
$\mathrm{Vdeg}(M)=\mathrm{deg}_{X_d}(M)$). Hence $M=M'$. Thus
the least monomials of two expansions of $g$ (with respect to the
ordering of Corollary \ref{well ordering}) are equal. Subtracting this monomial from two representations, and iterating the last procedure we deduce that  these two
expansions are the same and we are done (An argument similar to the last part  works for the initial of the induction $d=1$).
 $\hfill \Box $
\begin{remark} {\rm  Given an SKP $[U_{i,j},\beta_{i,j}]_{i=0..d,j=1..\alpha_i},$ and an
element
 $f\in \ke^{(d)}$ in order to obtain its  $(U)_{\alpha}-adic$ expansion,  we can use the algorithm of
  Proposition \ref{alg} for the acceptable vectors $\alpha '=(1,\ldots,1)$ and $\alpha ''=\alpha.$ Notice that in this case the
  $(U)_{\alpha '}-adic$ expansion of every  element $f\in \ke^{(d)}$ is
  itself.
 }
\end{remark}

 We also use the notation of $(\alpha ')-adic$ expansion. When there is
no stress on the specific  acceptable vector $\alpha '$ or it is understood,
we will talk about $U_d-adic$ or $adic$ expansion.

}

\section{Valuations associated to SKP's}{
 In this section we show that  to any SKP one can associate a valuation $\nu$ of the field $k((X_0,\ldots,X_d))$ centered on the
 ring $k[[X_0,\ldots,X_d]].$
\begin{definition}\label{valuation}{\rm
Let $[U_{i,j},\beta_{i,j}]$ be an SKP. For an acceptable
vector $\alpha ',$  we   define a map
$$\nu_{\alpha '}:\ \ke^{(d)}\backslash\{ 0\}\to \Phi$$ by:

\begin{itemize}
\item If $M$ is  any monomial $M(U)$ with $(U)_{\alpha '}-adic$ expansion
$M=c.U^{\mathbf{p}},$ where $c\in k$ then
$$\nu_{\alpha '}(M)=\sum_{i=0}^d\sum_{j=0}^{\alpha '_i}p_{i,j}\beta_{i,j}.$$
\item If $f\in \ke^{(d)}$ has the $ (U)_{\alpha '}-adic $ expansion
$f=\sum_{I(J)}c_{I(J)}U^{I(J)}$ then
$$\nu_{\alpha '}(f)=\mathrm{min}_{ I(J) }\{\nu_{\alpha '}(U^{I(J)})\}.$$
\end{itemize}
For any SKP, we denote the mapping of the definition  above
  by
$\nu_{\alpha}=\mathrm{val}[U_{i,j},\beta_{i,j}].$ We  will
see that this mapping is a valuation  (Theorem \ref{valdef}).

}
\end{definition}
\begin{definition}\label{ini-defi}{\rm
Let $[U_{i,j},\beta_{i,j}]$ be an SKP and   $f\in
\ke^{(d)}$ an arbitrary element and let $(\alpha ')$ be an
acceptable vector for this SKP. The initial form of $f$ with respect to
$\nu_{\alpha '}$ is defined as:
$$\mathrm{in }_{\nu_{\alpha '}}(f)=\sum_{I(J^{0})} c_{I(J^0)}U^{I(J^0)},$$
where $f=\sum_{I(J)}c_{I(J)}U^{I(J)}$ is the $(U)_{\alpha '}-adic$
expansion of $f$ and $I(J^0)$ ranges over those indices with
minimal $\nu_{\alpha '}-$value. }
\end{definition}

\begin{definition}\label{vector power}{\rm
Let $[U_{i,j},\beta_{i,j}]_{i=0..d,j=1..\alpha_i}$ be an
SKP and consider the power series ring $k_{((\alpha,d))}$. For any
monomial $M(U)=U^{\mathbf{a}}\in k_{((\alpha,d))}$ we define the vectors of the powers
$$\mathrm{VP}(M(U))=(\mathbf{a}_{d,\alpha_d},\mathbf{a}_{d-1,\alpha_{d-1}},\ldots,\mathbf{a}_{0,\alpha_0})\in \mathbb{N}^{d+1}.$$ }
\end{definition}

\begin{lemma}\label{max-degree}
Fix an {\rm SKP} and suppose that $\alpha '$ is an acceptable vector for
this {\rm  SKP.} Let $f\in \ke^{(d)}$ and suppose
$\mathrm{in}_{\nu_{\alpha '}}(f)=\sum_{I(J)}c_{I(J)}U^{I(J)}$. Then the
vectors of the powers $\mathrm{VP}(M)$  of the monomials $M$ of
$\mathrm{in}_{\nu_{\alpha '}}(f)$ are all different.
\end{lemma}

\noindent {\bf Proof. } Let $cU^{I(J)}$ and $c'U^{I'(J')}$ be two
monomials of $\mathrm{in}_{\nu_{\alpha '}}(f)$ with equal vectors of the powers. We show that for any $j=1,\ldots,\alpha '_d$ the powers of
the $U_{d,j}$ in the two monomials are the same. Indeed, let $j'<\alpha '_d$ be
the greatest index such that $I(J)_{d,j'}\neq I'(J')_{d,j'}$, note that this maximum index exists. We
assume $I(J)_{d,j'}>I'(J')_{d,j'}.$ By equating the
$\nu_{\alpha '}-$values of the two monomials
$$(I(J)_{d,j'}-I'(J')_{d,j'})\beta_{d,j'}=\sum_{(i'',j'')<_{lex}(d,j')}-(I(J)_{i'',j''}-I'(J')_{i'',j''})\beta_{i'',j''}\in(\beta_{i'',j''})_{(i'',j'')<_{lex}(d,j')},$$
which  is clearly a contradiction, because
$0<I(J)_{d,j'}-I'(J')_{d,j'}<n_{d,j'}.$ Continuing similar
argument for $i<d$, we deduce that the two monomials are the
same.$\hfill \Box$

\begin{corollary}\label{one initial}
Fix an {\rm SKP} and suppose $U_{i,\alpha_i}=0$, for $i=1..d$. For an arbitrary $0\neq f\in \ke^{(d)}$ the initial $\mathrm{in}_{\nu_{\alpha}}(f)$ consists of just one monomial of adic form.
\end{corollary}

\begin{lemma} \label{tec}
 Fix an { \rm SKP}
 $[U_{i,j},\beta_{i,j}]$ and let
$\alpha '$ be an acceptable vector. For any arbitrary  monomial $M(U)\in
k_{((\alpha ',d))},$ where $M=c.U^{\mathbf{ a}},$ we have:
\begin{itemize}
\item[$(i)$] The initial form of $M$ in its $(U)_{\alpha '}-adic$
expansion is just one monomial $M'=c'U^{\mathbf{ a'}}.$ In other
words, we have $\mathrm{in }_{\nu_{\alpha '}}(M)=M'.$
\item[$(ii)$]
We have
  $\mathbf{a'}_{d,\alpha '_d}=\mathbf{a}_{d,\alpha '_d}.$
\item[$(iii)$] For any two monomials $M$ and $M'$ of the power series ring $k_{((\alpha ',d))}$ with equal $\nu_{\alpha '}-$values, if
$\mathrm{VP}(M)<_{lex}\mathrm{VP}(M')$ then
$\mathrm{VP}(\mathrm{in}_{\nu_{\alpha '}}(M))<_{lex}\mathrm{VP}(\mathrm{in}_{\nu_{\alpha '}}(M'))$.

\end{itemize}
\end{lemma}

\noindent {\bf Proof. } For the first claim,  let $U_{i,j}$ be a
factor of  $M$ with power greater than $n_{i,j}$. Replace
$U_{i,j}^{n_{i,j}}$ by its expression from the algorithm for getting
$adic$ expansion. The claim is that after  one such replacement
there exists just one monomial with minimal $\nu_{\alpha '}-$value. We
prove the claim for the replacements of the first type of
algorithm for getting $adic$ expansion. For the second type the
argument is similar. After a replacement of type one we get two
monomials with different $\nu_{\alpha '}-$values:
\[\begin{array}{ll }
 M& =\frac{M}{U_{i,j}^{n_{i,j}}}(U_{i,j+1} +\theta_{i,j}U^{m^{(i,j)}})\\
&\begin{array}{ccc}
=\underbrace{c\frac{U^{\mathbf{a}}.U_{i,j+1}}{U_{i,j}^{n_{i,j}}}}&+&\underbrace{c\theta_{i,j}\frac{U^{\mathbf{a}}U^{m^{(i,j)}}}
{U_{i,j}^{n_{i,j}}}}.\\
 \mbox{$M_2$}& &\mbox{$M_1$}
\end{array}\\
\end{array}\]
Then $\nu_{\alpha '}(M_2)> \nu_{\alpha '}(M_1)=\nu_{\alpha '}(M).$ Therefore, we
have $\mathrm{in }_{\nu_{\alpha '}}(M)=\mathrm{in}_{\nu_{\alpha '}}(M_1).$   We do the same for $M_1.$ Finally we get a  monomial $M'$ whose
 $adic$ expansion is  itself, this proves $(i).$

 For the second part we notice the that the proof of the first part shows  the following  general fact:\\
 \noindent
 For the monomial $M(U)$ a replacement on $U_{i,j}^{n_{i,j}}$ cannot affect the power of $U_{i',j'}$, for
 $(i',j')>_{lex.}(i,j)$,
 of the unique monomial with minimal value of the expansion generated after replacement.

 For the last part, suppose $M=U^{\mathbf{a}}$ and $M'=U^{\mathbf{a}'}$. Let $d'\leq d$ be the first index such that
 $\mathbf{a}_{d',\alpha '_{d'}}<\mathbf{a}'_{d',\alpha '_{d'}}$. Then by Lemma \ref{max-degree}, we have
 $\mathbf{a}_{i,j}=\mathbf{a}'_{i,j}$ for $i=d'+1,\ldots,d$ and $j=1,\ldots,\alpha_i$.
 Thus the algorithm for  getting $adic$ expansion for these two monomials for such $i$ and $j$ can be chosen the same. Hence, without loss of
 generality we can assume that $\mathbf{a}_{i,j}<n_{i,j}$ and $\mathbf{a}'_{i,j}<n_{i,j}$, for $i=d'+1,\ldots,d$ and $j<\alpha_i$. Then because
 $\mathbf{a}_{d',\alpha '_{d'}}<\mathbf{a}'_{d',\alpha '_{d'}}$, by part $(ii)$ we are done.    $\hfill\Box$

\begin{theorem} \label{valdef}
Given any {\rm SKP} $[U_{i,j},\beta_{i,j}],$ for any
acceptable vector  $\alpha ',$ the mapping\\ $\nu_{\alpha '}:\
\ke^{(d)}\backslash \{0\}\to \Phi$ extends in an obvious way  to
a $k-$valuation of the field $k((X_0,\ldots,X_d)).$ Moreover, for
any two acceptable vectors $\alpha '$ and $\alpha ''$ such that $\alpha '\leq \alpha ''$
and for any $f\in \ke^{(d)}$ we have
$\nu_{\alpha '}(f)\leq\nu_{\alpha ''}(f).$

 \end{theorem}

\noindent
 {\bf Proof. } The extension to the field $k((X_0,\ldots,X_d))$ is a trivial
 task. We need only to prove  that given any $f,g\in \ke^{(d)}\backslash\{ 0\}$
 we have $\nu_{\alpha '} (f+g)\geq
 \mathrm{min}\{\nu_{\alpha '}(f),\nu_{\alpha '}(g)\}$ and
 $\nu_{\alpha '}(f.g)=\nu_{\alpha '}(f)+\nu_{\alpha '}(g).$ The first one
 is a direct consequence of the definition and the uniqueness of
 the $adic$ expansions. For the second equality, let
 $\mathrm{in}(f)=\sum_{I(J)}c_{I(J)}U^{I(J)}$ and
 $\mathrm{in}(g)=\sum_{I'(J')}c'_{I'(J')}U^{I'(J')}.$ Let
 $M= c_{I(J_0)}U^{I(J_0)}$ (respectively  $M'= c_{I'(J'_0)}U^{I'(J'_0)}$) be the unique (Lemma \ref{max-degree})  monomial of
 the expansion of $\mathrm{in}(f)$ (respectively $\mathrm{in}(g)$) with minimal vector of powers,  with
 respect to the $lex.$ order.
  Then by Lemma \ref{tec}, $(iii)$,
 we see that $\mathrm{in}(M.M')=M'' $ is the unique  monomial of $\mathrm{in}(f.g)$ with
 minimal
  vector of the  powers.
      But $\nu_{\alpha '}(M'')=
 \nu_{\alpha '}(M)+\nu_{\alpha '}(M')=\nu_{\alpha '}(f)+\nu_{\alpha '}(g),$
  by the definition of the mapping $\nu_{\alpha '}$ we have
 $\nu_{\alpha '}(M'')=\nu_{\alpha '}(f.g).$ For the last part, we note that in
 the algorithm for getting $\alpha ''-adic$ expansion of an element from
 its $\alpha '-adic$ expansion at every step in the substitution we replace a
 monomial with two new monomials with values  equal to
or greater than the original monomial.
 %The last part is an easy consequence of the algorithm of getting $(\alpha '')-adic$
 %expansion of $f$ from its $(\alpha ')-adic$ expansion.
 $\hfill \Box$
\begin{corollary}
Given an {\rm SKP} $[U_{i,j},\beta_{i,j}]_{i=0..d,j=1..\alpha_i}$, all the $U_{i,j}$'s  are irreducible elements
of the power series ring $\ke^{(i-1)}[X_i].$
\end{corollary}
{\bf Proof. } We prove the claim for $U_{d,j}$. Consider the
vector $(\alpha '),$ defined by $\alpha '_i=\alpha_{i}$, for $0\leq i< d$, and
$\alpha '_d=j$. This is an acceptable vector. In this proof all the
$adic$ expansions are $(U)_{\alpha '}-adic$ expansions. We give a proof by contradiction. Assume that $U_{d,j}$ is reducible and $U_{d,j}=f.g,$ for some non-unit elements $f,g\in
\ke^{(d-1)}[X_d].$
    As the $\alpha '-adic$ expansion of $U_{d,j}$ is  itself, we have $\mathrm{in}(U_{d,j})=U_{d,j}.$
We can compute this  initial in the other way, using  initials of
$f$ and $g.$ This gives us $U_{d,j}=
\mathrm{in}(\mathrm{in}(f).\mathrm{in}(g)).$

On the other hand,
$\beta_{d,j}=\nu_{\alpha '}(U_{d,j})=\nu_{\alpha '}(f)+\nu_{\alpha '}(g)$.  Thus the
monomials of   $\mathrm{in}(f)$ and $\mathrm{in}(g))$ do not
have a factor $U_{d,j}$. By Lemma \ref{tec}, $(ii)$, this shows
that the monomials $\mathrm{in}(\mathrm{in}(f).\mathrm{in}(g))$
do not have a factor $U_{d,j}$, which is a contradiction. $\hfill \Box$

\begin{remark}\label{diff-SKP}{\rm
One should note that in the definition of the SKP's for the ring
$k[[X_0,\ldots,X_d]]$ the ordering of the variables plays an
important role. In other words, changing the coordinates of the
rings (even with a permutation) may change totally the system of
SKP's associated to the valuation, or even they may not exist. This phenomenon can be seen even
in dimension two; For example consider the valuation $\nu$
centered
 on the ring $k[X_0,X_1]$ defined by the SKP,
$[(U_{0,1},U_{1,1},U_{1,2},U_{1,3}),(2,3,9,10)].$ Where, we have
$U_{1,2}=U_{1,1}^2-U_{0,1}^3,\ U_{1,3}=U_{1,2}-U_{0,1}^3U_{1,1}. $
Note that the last two equations are given to us (up to the
knowledge of the corresponding $\theta$'s ) as soon as the
sequence of $\beta$'s $(2,3,9)$ is known. Now, changing the order
of the coordinates, we consider the same ring as $k[Y_0,Y_1]$ with
$Y_0=X_1,Y_1=X_0.$ The same valuation is given by the following
SKP's in the new coordinate
$\nu=\mathrm{val}[(V_{0,1},V_{1,1},V_{1,2}, V_{1,3}),(3,2,9,10)].$
Where the SKP's are as follows.
$$V_{1,2}=V_{1,1}^3-V_{0,1}^2,\ V_{1,3}=V_{1,2}+V_{0,1}^3.$$
The relation between two SKP's is as follows:
$$V_{0,1}=U_{1,1},\ V_{1,1}=U_{0,1},\ V_{1,2}=-U_{1,2}.$$
For   $V_{1,3}$ we have:
$$V_{1,3}=V_{1,2}+V_{0,1}^3=-U_{1,2}+U_{1,1}^3=-U_{1,2}+(U_{0,1}^3+U_{1,2})U_{1,1}=-U_{1,3}+U_{1,1}U_{1,2}.$$
As this example shows the explicit relation between the $U$'s and
$V$'s is not,  in general, trivial.}
\end{remark}

}

\section{Euclidean expansion and other properties of SKP's}{
In this section we give another expansion in the  ring
$k_{(d-1)}[X_d]$, associated to a SKP of the power series ring
$k[[X_0,\ldots,X_d]]$ ($k_{(i)}:=k((X_0,\ldots,X_i))$).   We show that the valuation $\nu$ associated to
this SKP, can be defined using this new expansion, plus the
knowledge of the valuation $\nu$ on the field $k_{(d-1)}$.
Moreover, we show that the Euclidean expansion can be obtained
directly from the $adic$ expansion. This is interesting in practice,
because $adic$ expansion is defined only with substitutions while
Euclidean expansion is defined using divisions.

\begin{definition}\label{Euclidean expansion} {\rm ({\bf Euclidean
expansion}) Fix an SKP
$[U_{i,j},\beta_{i,j}]_{i=0..d,j=1..\alpha_i}.$ For
any $j=1,\ldots,\alpha_d$ we define the  acceptable vector
$\alpha^{(j)}=(\alpha_0,\ldots,\alpha_{d-1},j)$. Let $f\in k_{(d-1)}[X_d]$, and
consider the expansion $f=\sum_{J}c_{J}U_{d}^{J}\in k_{(d-1)}[U_{d}]$ such that $0\leq J_{j'}<n_{d,j'}$ for any $0\leq j'<j$.
This is called the $j$th Euclidean expansion of
$f$.}\end{definition}

\begin{proposition}{\rm\bf (Algorithm for getting Euclidean
expansion)}  With the notations of Definition \ref{Euclidean
expansion}, do the following:

Consider the greatest index $j_0$ such that
$\mathrm{deg}_{X_d}(f)>d_{d,j_0}$. Divide $f$ by $U_{d,j_0}$ in
the ring $k_{(d-1)}[X_d]$ to obtain $f=qU_{d,j_0}+r,$ where
$q,r\in k_{(d-1)}[X_d]$ and $\mathrm{deg}_{X_d}(r)<d_{d,j_0}$.
Iterate the same process for $q$ as far as possible to obtain
$f=\sum_{t}f_{t}U_{d,j_0}^{t},$ where
$\mathrm{deg}_{X_d}(f_{t})<d_{d,j_0}$. Iterate the same procedure
for each of the $f_{t}$'s and the greatest index $j'$, $j'<j_0$,
such that $d_{d,j'}<d_{d,j_0}$. Continue as far as possible.
This process terminates after finitely many steps. The
resulting expansion is equal to the $j${\rm th} Euclidean
expansion of $f$. Moreover, the Euclidean expansion is unique.
\end{proposition}
\noindent {\bf Proof. }  As the $U_{d}$'s which appear in the
process are among the elements of the finite set $\{U_{d,j'}/ \
n_{d,j'}\neq 1,\ \mathrm{and\ } \mathrm{deg}_{X_d}(f)>d_{d,j'}\}$,
the process stops after finitely many steps. We show that the
resulting expansion is the $j$th Euclidean expansion of $f$. Let
$U_{d}^{J}$ be a monomial generated in the algorithm above. It is
sufficient to show that this monomial is of Euclidean form.
Indeed, let $j'$ be the greatest index less than $j$ such that
$J_{j'}\geq n_{d,j}$. This means that
$\mathrm{deg}_{X_d}(U_{d,1}^{J_1}\cdots U_{d,j'}^{J_{j'}})\geq
d_{d,j'+1}$ and we must divide it (in the monomial in the
procedure above) by $U_{d,j'+1}$, which is a contradiction.

The uniqueness of Euclidean expansion comes from the fact that (by
Lemma \ref{deg-rev}) the $\mathrm{deg}_{X_d}(U_{d}^{J})$ of a
monomial of Euclidean form determines the vector $J$. Therefore
there is a unique vector $J_0$ such that
$\mathrm{deg}_{X_d}(U_{d}^{J_0})=\mathrm{deg}_{X_d}(f)$. This
monomial (plus its coefficient)is common in all the possible
Euclidean expansions of $f$. Subtracting this monomial from $f$,
by induction on the degree of $f$ we are done. $\hfill\Box$

\begin{lemma}
Fix an {\rm SKP }$ [U_{i,j},\beta_{i,j}]_{i=0..d,j=1..\alpha_i}$ and let  $f\in k[[X_0,\ldots,X_d]]$. The $j${\rm
th} Euclidean expansion of $f$ can be obtained using the
$(\alpha^{(j)})-adic$ expansion of it as follows. Let
$f=\sum_{I(J)}c_{I(J)}U^{I(J)}$ be the $(\alpha^{(j)})-adic$ expansion
of $f$. Then the Euclidean expansion of $f$ is equal to
$$\sum_{J'}(\sum_{I(J),I(J)_{d}=J'}c_{I(J)}\frac{U^{I(J)}}{U_{d}^{J'}})U_{d}^{J'}.$$
\end{lemma}
\noindent {\bf Proof. } It is clear that the above expansion
satisfies all the properties of the $j$th Euclidean expansion of
$f$. Thus, by uniqueness, it is the Euclidean expansion of $f$.
$\hfill\Box$
%%%%%%%%%%%%%%%%%%%%%%%%%%%%%%%%%%%%%%%%%%%%%%%%%%%%%%%%%%%%%%%%%%%%%%%%%%%%%%%%%%%%%%%%%%%%%%%%%%%%%%%%%%%%%%%%%%%%%%%%%%%%%%
\begin{remark}\label{extended Euclidean expansion}
{\rm Using  the  above lemma, we  extend the notion of Euclidean expansion to the power series ring $\ke^{(d)}$.  An expansion of  $f\in \ke^{(d)}$  of the form $f=\sum_{J}c_JU_d^J\in \ke^{(d-1)}[[U_d]]$ which
satisfies the conditions  of Definition \ref{Euclidean expansion} is called the Euclidean expansion of $f$. The above lemma
shows that such an expansion can be obtained using $adic$ expansion of $f$. An argument, similar to the proof of Proposition \ref{alg}, shows that this expansion is unique. }
\end{remark}

\begin{proposition}
Fix an {\rm SKP }$ [U_{i,j},\beta_{i,j}]_{i=0..d,j=1..\alpha_i}$ and let $\nu$ be the $k-$valuation of the field
$k_{(d)}$ associated to it. Set
$\overline{\nu}=\nu\mid_{k_{(d-1)}}$. The valuation $\nu$ (as a
valuation of the field $k_{(d-1)}(X_d)$) can be defined using the
data
$[\overline{\nu},(U_{d,j})_{j=1}^{\alpha_d},(\beta_{d,j})_{j=1}^{\alpha_d}]$
as follows. For any $f\in \ke^{(d-1)}[X_d]$ let
$f=\sum_{J}f_{J}U_{d}^{J}$ be its $\alpha_d${\rm th} Euclidean
expansion then
$$\nu(f)=\mathrm{min}_{J}\{\overline{\nu}(f_{J})+\beta_{d}.J\}.$$
\end{proposition}

\noindent {\bf Proof. } The lemma above shows that the equation of
the proposition is just another way of writing $\nu(f)$, which is
originally the minimum of the values of the  monomials in the  $adic$
expansion of $f$. $\hfill\Box$\\

\begin{remark}\label{Favre Jonsson def}{\rm With the notations of the proposition above, write $f=\sum_{t}f_tU_{d,j}^t,$
 with  $\mathrm{deg}_{X_d}(f_t)<d_{d,j}.$ Then with a similar argument we
have $$\nu(f)=\min_{t}\{\nu(f_t)+t\beta_{d,j}\}.$$
 }
\end{remark}

\begin{definition}\label{delta-Euclid}{\rm
Fix an {\rm SKP} $[U_{i,j},\beta_{i,j}]_{i=0..d,j=1..\alpha_i}.$ We consider the set of acceptable vectors
$\alpha^{(j)}=(\alpha_0,\ldots,\alpha_{d-1},j),$ for $j=1,\ldots, \alpha_d.$

For any $f\in k_{(d-1)}[X_d],$  and any $\alpha^{(j)}$ we define
$$\delta _{\alpha^{(j)}}(f)=\mathrm{max}\{\ell:\ \ell \mathrm{ \ is\
power\ of\ } U_{d,j} \mathrm{ \ in\ the\ monomials\ of\
in}_{\nu^{\alpha^{(j)}}}(f)\}.$$}
\end{definition}

\begin{remark}
{\rm
Let  $u\in k^{(d-1)}$ and $f\in k_{(d-1)}[X_d]$  then $\delta _{\alpha^{(j)}}(f)=\delta _{\alpha^{(j)}}(uf)$.
}
\end{remark}

\begin{lemma}\label{delta}
For any $f,g\in k_{(d-1)}[X_d]$   we have
$$\delta_{\alpha^{(j)}}(f.g)=\delta_{\alpha^{(j)}}(f)+\delta_{\alpha^{(j)}}(g).$$
\end{lemma}

\noindent {\bf Proof. } First we find $u,v\in k^{(d-1)}$ such that $uf,vg\in k^{(d-1)}[X_d]$, this is always possible. Then by last remark it suffices to prove the lemma for $uf$ and $vg$, i.e., we can assume $f,g\in k^{(d-1)}[X_d]$.  Lemma \ref{max-degree} shows that  there are
unique monomials $f_{J}U_d^{J}$ and $g_{J'}U_d^{J'}$ of
$\mathrm{in}(f)$ and $\mathrm{in}(g)$ (respectively) that have
maximal $U_{d,j}$ power. Write Euclidean expansion of
$\mathrm{in}(f.g)$ using the product
$\mathrm{in}(f).\mathrm{in}(g)$ and algorithm for getting $adic$
expansion. We see $\mathrm{in}(f).\mathrm{in}(g)$ has a unique
monomial with $U_{d,j}-$degree equal
$\delta_{\alpha^{(j)}}(f)+\delta_{\alpha^{(j)}}(g),$ i.e.,
$f_{J}g_{J'}U^{J}U^{J'}.$ Now, Lemma \ref{tec}, $(ii)$, shows that
after getting $adic$ expansion from this product the
$U_{d,j}-$powers of the monomials do not change which proves the
equality. $\hfill \Box$\\

The following lemma is an adaptation  of the results of \cite{FJ} in
our situation.

\begin{lemma}\label{old-results}
Fix an {\rm SKP} $[U_{i,j},\beta_{i,j}],$ and let $\alpha^{(j)}$ be defined as in Definition \ref{delta-Euclid} then
\begin{itemize}
\item[($i$)] For $f\in k_{(d-1)}[X_d],$ we have $\delta_{\alpha^{(j)}}(f)=0$ iff $\mathrm{in}_{\nu_{\alpha^{(j)}}}(f)$ is
a unit in $\mathrm{gr}_{\nu_{\alpha^{(j)}}}k_{(d-1)}[X_d].$
\item[(ii)] If $f,g\in k_{(d-1)}[X_d]$ then there exist $Q,R\in k_{(d-1)}[X_d]$
such that $\mathrm{in}_{\nu_{\alpha^{(j)}}} (f)=\mathrm{in}_{\nu_{\alpha^{(j)}}} (Qg+R)$ in $\mathrm{gr}_{\nu_{\alpha^{(j)}}}k_{(d-1)}[X_d]$
and $\delta_{\alpha^{(j)}}(R)<\delta_{\alpha^{(j)}}(g).$
\item[($iii$)] The polynomials  $\mathrm{in}_{\nu_{\alpha^{(j)}}}(U_{d,\alpha^{(j)}_j}),\mathrm{in}_{\nu_{\alpha^{(j)}}}(U_{d,\alpha^{(j)}_{j+1}})$ are irreducible in $\mathrm{gr}_{\nu_{\alpha^{(j)}}}k_{(d-1)}[X_d].$
\item[($iv$)] If $j'<j$ then $\mathrm{in}_{\nu_{\alpha^{(j)}}}(U_{d,j'})$ is a unit in $\mathrm{gr}_{\nu_{\alpha^{(j)}}}k_{(d-1)}[X_d].$
\item[($v$)] If $f=\sum_{t}f_tU_{d,j}^t,$
 with  $\mathrm{deg}_{X_d}(f_t)<d_{d,j}$ and $\delta_{\alpha^{(j)}}(f)<n_{d,j}$ then $\mathrm{in}_{\nu_{\alpha^{(j)}}}(f)=\mathrm{in}_{\nu_{\alpha^{(j)}}}(f_tU_{d,j}^t)$ in
$\mathrm{gr}_{\nu_{\alpha^{(j)}}}k_{(d-1)}[X_d]$ for some $t<n_{d,j}.$
 \end{itemize}
\end{lemma}

\noindent {\bf Proof. } Throughout the proof we fix the expansion
$f=\sum_{t}f_tU_{d,j}^t,$
 with  $\mathrm{deg}_{X_d}(f_t)<d_{d,j}.$\par\noindent
 $(i).$  If $\delta_{\alpha^{(j)}}(f)=0$ then $\mathrm{in}_{\nu_{\alpha^{(j)}}}(f)=\mathrm{in}_{\nu_{\alpha^{(j)}}}(f_0)$ in
$\mathrm{gr}_{\nu_{\alpha^{(j)}}}k_{(d-1)}[X_d]$. As $U_{d,j}$ is
irreducible and $\mathrm{deg}_{X_d}(f_0)<d_{d,j}$ the polynomial
$U_{d,j}$ is prime with $f_0$. Hence  we can find $A,B\in
k_{(d-1)}[X_d],$
$\mathrm{deg}_{X_d}(A),\mathrm{deg}_{X_d}(B)<d_{d,j}$ so that
$Af_0=1-BU_{d,j}.$ Then $\nu_{\alpha^{(j)}}(Af_0)=\nu_{\alpha^{(j)}}(1)<
\nu_{\alpha^{(j)}}(BU_{d,j}).$ Therefore, $\mathrm{in}_{\nu_{\alpha^{(j)}}}(Af_0)=1$ in
$\mathrm{gr}_{\nu_{\alpha^{(j)}}}k_{(d-1)}[X_d]$. So $\mathrm{in}_{\nu_{\alpha^{(j)}}}(f_0)$ and hence
$\mathrm{in}_{\nu_{\alpha^{(j)}}}(f)$ is a unit in $\mathrm{gr}_{\nu_{\alpha^{(j)}}}k_{(d-1)}[X_d]$.
Conversely, if $\mathrm{in}_{\nu_{\alpha^{(j)}}}(f)$ is unit, say $\mathrm{in}_{\nu_{\alpha^{(j)}}}(Af)=1$ in
$\mathrm{gr}_{\nu_{\alpha^{(j)}}}k_{(d-1)}[X_d]$ for some $A\in
k_{(d-1)}[X_d]$ then
$\delta_{\alpha^{(j)}}(f)+\delta_{\alpha^{(j)}}(A)=\delta_{\alpha^{(j)}}(1)=0$
 so $\delta_{\alpha^{(j)}}(f)=0$.\par\noindent
$(ii).$ Write $g=\sum_{t}g_tU_{d,j}^{t}$. It suffices  to prove
the claim when $g_t=0$ for $t>M:=\delta_{\alpha^{(j)}}(g)$ and using
$(i)$ we may assume $g_{M}=1$. As
$\mathrm{deg}_{X_d}(g_t)<d_{d,j}$ for $t\leq M$ we have
$\mathrm{deg}_{X_d}(g)=Md_{d,j}$. Euclidean division in
$k_{(d-1)}[X_d]$ yields $Q,R^1\in k_{(d-1)}[X_d]$ with
$\mathrm{deg}_{X_d}(R^1)<\mathrm{deg}_{X_d}(g)$ so that
$f=Qg+R^1$. Write $R^1=\sum_{t}R_{t}U_{d,j}^t$ and set
$N:=\delta_{\alpha^{(j)}}(R^1),$ $R:=\sum_{t\leq N}R_{t}U_{d,j}^t$.
Then $\mathrm{in}_{\nu_{\alpha^{(j)}}}(f)=\mathrm{in}_{\nu_{\alpha^{(j)}}}(Qg+R)$ in $\mathrm{gr}_{\nu_{\alpha^{(j)}}}k_{(d-1)}[X_d]$ and
$$\mathrm{deg}_{X_d}(R)=\mathrm{deg}_{X_d}(R_N)+Nd_{d,j}<Md_{d,j}=\mathrm{deg}_{X_d}(f).$$
Hence $N<M$ and we are done.\par\noindent
$(iii).$ We have
$\delta_{\alpha^{(j)}}(U_{d,j})=1$ so if $\mathrm{in}_{\nu_{\alpha^{(j)}}}(U_{d,j})=fg$ in
$\mathrm{gr}_{\nu_{\alpha^{(j)}}}k_{(d-1)}[X_d]$ then
$\delta_{\nu_{\alpha^{(j)}}}(f)=0$ or $\delta_{\alpha^{(j)}}(g)=0.$ Hence by
$(i)$, $\mathrm{in}_{\nu_{\alpha^{(j)}}}(f)$ or $\mathrm{in}_{\nu_{\alpha^{(j)}}}(g)$ is a unit in
$\mathrm{gr}_{\nu_{\alpha^{(j)}}}k_{(d-1)}[X_d]$. \par\noindent For
$U_{d,j+1},$ we have
$U_{d,j+1}=U_{d,j}^{n_{d,j}}-\theta_{d,j}U^{m^{(d,j)}}$. Let
$\mathrm{in}_{\nu_{\alpha^{(j)}}}(U_{d,j+1})=\mathrm{in}_{\nu_{\alpha^{(j)}}}(fg)$ in $\mathrm{gr}_{\nu_{\alpha^{(j)}}}k_{(d-1)}[X_d]$ with
$0<\delta_{\alpha^{(j)}}(f),\delta_{\alpha^{(j)}}<n_{d,j}$. By $(v)$, we can
write $f=f_tU_{d,j}^t,\ g=g_{t'}U_{d,j}^{t'}.$ Then
$U_{d,j+1}=f_tg_{t'}U_{d,j}^{n_{d,j}}$ so
$(1-f_tg_{t'})U_{d,j}^{n_{d,j}}=\theta_{d,j}U^{m^{(d,j)}}$. As
$U_{d,j}$ is irreducible and $U^{m^{(d,j)}}$ a unit, we have
$\mathrm{in}_{\nu_{\alpha^{(j)}}}(f_tg_{t'})=1$ in $\mathrm{gr}_{\nu_{\alpha^{(j)}}}k_{(d-1)}[X_d]$. But
then $\mathrm{in}_{\nu_{\alpha^{(j)}}}(U^{m^{(d,j)}})=0$ in
$\mathrm{gr}_{\nu_{\alpha^{(j)}}}k_{(d-1)}[X_d]$ which is absurd. So we
can assume $\delta_{\alpha^{(j)}}(f)=n_{d,j}$ and $\delta_{\alpha^{(j)}}=0.$
Hence $g$ is a unit.\par\noindent
 $(iv).$ By $(i)$ it suffices to show that $\delta_{\alpha^{(j)}}(U_{d,j'})=0.$
If $d_{d,j'}<d_{d,j}$ then this is
 obvious. If $d_{d,j'}=d_{d,j}$ then $U_{d,j'}=(U_{d,j'}-U_{d,j})+U_{d,j}$
 where $\mathrm{deg}_{X_d}(U_{d,j'}-U_{d,j})<d_{d,j}.$ Now $\nu_{\alpha^{(j)}}(U_{d,j'})=\beta_{d,j'}<\beta_{d,j}=\nu_{\alpha^{(j)}}(U_{d,j}),$
 so $\nu_{\alpha^{(j)}}(U_{d,j'}-U_{d,j})<\nu_{\alpha^{(j)}}(U_{d,j})$ and
$\delta_{\alpha^{(j)}}(U_j')=0.$\par\noindent
$(v).$  Suppose
$\nu_{\alpha^{(j)}}(f_tU_{d,j}^t)=\nu_{\alpha^{(j)}}(f_{t'}U_{d,j}^{t'}),$
where $t\leq t'<n_{d,j}.$ Then $(t'-t)\beta_{d,j}=
\nu_{\alpha^{(j-1)}}(f_t)-\nu_{\alpha^{(j-1)}}(f_{t'}).$ Hence $n_{d,j}\mid t'-t$ thus $t'=t.$ $\hfill\Box$\\

\begin{proposition}\label{proposition:euclidean domain}
The graded algebra  $\mathrm{gr}_{\nu_{\alpha}}k_{(d-1)}[X_d]$ is a
Euclidean domain.
\end{proposition}
\noindent{\bf Proof. } The item $(ii)$  of Lemma \ref{old-results} proves the claim. $\hfill\Box$

\begin{theorem}\label{lemma representation graded}
Fix an {\rm SKP} $[U_{i,j},\beta_{i,j}]_{i=0..d,j=1..\alpha_i}$ and let $\nu$ be its associated valuation.
Consider $0\neq f\in k[[X_1,\ldots,X_d]].$ Then initial form of  $f$ has a unique
decomposition of the form:
\begin{itemize}
\item[$(i)$] If $U_{d,\alpha_d}\neq 0,\ n_{d,\alpha_d}=\infty$ then
 $$f=\tilde{f}U_{d}^{J},\ \mathrm{in}\
\mathrm{gr}_{\nu}k_{(d-1)}[X_d],
$$
where $\tilde{f}\in k_{(d-1)}$ and $0\leq J_j<n_{d,j},$ for $1\leq
j<\alpha_d.$
\item[$(ii)$] If $U_{d,\alpha_d}\neq 0,\ n_{d,\alpha_d}\neq \infty$
then $$f=p(T)U_{d}^{\hat{J}},\ \mathrm{in}\
\mathrm{gr}_{\nu}k_{(d)},$$ where $p(T)\in k_{(d-1)}[T]$ and
$0\leq J_j<n_{d,j},$ for $1\leq j\leq \alpha_d,$ and
$T=U_{d,\alpha_d}^{n_{d,\alpha_d}}U^{-m^{(d,\alpha_d)}}.$ Moreover, all the
coefficients of $p(T)$ have the same $\nu-$value.
\item[$(iii)$] If $U_{d,\alpha_d}=0$ then
$$f=\tilde{f}U_{d}^{J},\ \mathrm{in}\
\mathrm{gr}_{\nu}k_{(d-1)}[X_d],$$ where $0\leq J_j<n_{d,j},$ for
$1\leq j\leq \alpha_d$, and $J_j=0,$ except for a finite number of $j$.
\end{itemize}
\end{theorem}

\noindent
 {\bf Proof.} $(i).$  Suppose
$f=\sum_{J}f_{J}U_{d}^{J}$ is the Euclidean expansion of $f $ (Remark \ref{extended Euclidean expansion}),
 where $f_{J}\in k_{(d-1)},$ and $0\leq
J_j<n_{d,j}$ for $j<\alpha_d$ .  We claim that for any two $J$ and $J'$
we have $\nu(f_{J}U_{d}^{J})\neq \nu(f_{J'}U_{d}^{J'}).$ Indeed,
if we have equality, consider the greatest index $j_0$ such that
$J_{j_0}\neq J'_{j_0}.$ We have
$(J'_{j_0}-J_{j_0})\beta_{d,j_0}=\nu(f_J)-\nu(f_{J'})+\sum{j<j_0}(J_j-J'_j)\beta_{d,j}.$
Then as $j_0<\alpha_d$ (because $n_{d,\alpha_d}=\infty$), we have
$n_{d,j_0}\mid J_{j_0}-J'_{j_0}.$ Thus $J_{j_0}=J'_{j_0}$ which is
absurd.
\par\noindent
$(ii).$ We show that any monomial $f_{J}U_{d}^{J}$ of the
Euclidean expansion of $\mathrm{in}(f)$ is of the form
$\hat{f}_{J}T^{r_{\alpha_d}} U_{d}^{\hat{J}},$ $ \mathrm{in}\
\mathrm{gr}_{\nu}\ke^{(d)},$ for a fixed $\hat{J}$ such that $0\leq
\hat{J}_j<n_{d,j},$ for any $j$.\\
Fix $J,$ and make the Euclidean division $J_{\alpha_d}=r_{\alpha_d}
n_{d,\alpha_d}+\hat{J}_{\alpha_d},$ $0\leq \hat{J}_{\alpha_d}<n_{d,\alpha_d}$. And
write
$f_{J}U_{d}^{J}=\overline{f}_{J}U_{d,\alpha_d}^{\hat{J}_{\alpha_d}}T^{r_{\alpha_d}}U_{d}^{\mathbf{a}},$
with $\mathbf{a}:=J+r_{\alpha_d}.m_d^{(d,\alpha_d)}.$ As
$$U_{d,j}^{n_{d,j}}=\theta_{d,j}(U_{<d-1}^{m_{<d-1}^{(d,j)}})U_{d}^{m_d^{(d,j)}},
\ \mathrm{in}\ \mathrm{gr}_{\nu}k_{(d-1)}[X_d],$$ making the
Euclidean division $\mathbf{a}_{j}=r_{j}n_{d,j}+\hat{J}_{j},$
(with $0\leq \hat{J}_{j}<n_{d,j}$) for the greatest index $j$ such
that $\mathbf{a}_j\neq 0,$ we get $\prod_{j'\leq
j}U_{d,j'}^{\mathbf{a}_{j'}}=U_{d,j}^{\hat{J}_{j}}\prod_{j'<j}U_{d,j'}^{\mathbf{a}'_{j'}}$
with $\mathbf{a}'_{j'}\in\mathbb{N}.$ We finally get by induction,
a representation
$$f_JU_{d}^{J}=\hat{f}_JT^{r_{\alpha_d}}U_{d}^{\hat{J}},$$
where $0\leq \hat{J}_j<n_{d,j},$ for any $j.$ As $\nu(T)=0,$ with
an argument like in the final part of the case $(i)$  one can
argue to show   that $\hat{J}$ is the same for all $J$'s. Clearly,
the coefficients of $p$ has the same $\nu-$value.
\par\noindent
$(iii).$ This is similar to $(i).$$\hfill\Box$

\begin{corollary}\label{irreducible elements}
Let $\nu$ be a valuation as above.
\begin{itemize}
\item[$(i)$] If $U_{d,\alpha_d}\neq 0,\ n_{d,\alpha_d}\neq \infty$ the
only irreducible element of $\mathrm{gr}_{\nu}k_{(d-1)}[X_d]$ is
$U_{d,\alpha_d}$.
\item[$(ii)$] Assume that $k$ is an algebraically closed field and  $U_{d,\alpha_d}\neq 0,\ n_{d,\alpha_d}<\infty$ and that the following additional  condition is satisfied:
 for every two monomial  $U^{I},U^{J}\in k_{(d-1)}$ of $adic$ form we have  $U^{I}=U^{J}$ whenever  $\nu(U^{I})=\nu(U^{J})$.
 Then the
irreducible elements  of $\mathrm{gr}_{\nu}k_{(d-1)}[X_d]$ are of
the form $U_{d,\alpha_d}^{n_{d,\alpha_d}}-\theta U^{m^{(d,\alpha_d)}},$ for some
$\theta\in k.$
\item[$(iii)$] If $U_{d,\alpha_d}=0$ then
$\mathrm{gr}_{\nu}k_{(d-1)}[X_d]$ is a field.
\end{itemize}
\end{corollary}

\noindent
{\bf Proof.} $(i).$ Assume
$f\in\mathrm{gr}_{\nu}k_{(d-1)}[X_d]$ is irreducible. By $(i)$ of
the last theorem, $f=\tilde{f}U_{d}^{J}$. But $U_{d,j}$ is a unit
for $j<\alpha_d$ (by Lemma \ref{old-results}, $(iv)$), so $U_{d,\alpha_d}$
is the only irreducible element in
$\mathrm{gr}_{\nu}k_{(d-1)}[X_d]$ (Lemma \ref{old-results},
$(iii)$).
\par\noindent
$(ii).$ We use $(ii)$ of the last theorem. There we construct a
polynomial $p(T)\in k_{(d-1)}[T]$.  As we are working in the
graded ring, we can replace the coefficients of $p$ with their
initial, which by  assumption is a unique monomial $U^{I_0}\in
k_{(d-1)}$. Thus $P(T)=U^{I_0}p'(T),$ where $p'(T)\in k[T]$.
Factorize $p'(T)=\prod(T-\theta_{l})$,  modulo unit factors, we
hence get
$$f=U^{I_0}U_{d,\alpha_d}^{\hat{J}_{\alpha_d}-Ln_{d,\alpha_d}}\prod_{l}(U_{d,\alpha_d}^{n_{d,\alpha_d}}-\theta_{l}U^{m^{(d,\alpha_d)}}),$$
where $L=\mathrm{deg}(p)$. On the other hand Lemma
\ref{old-results}, $(iii)$, shows that all the elements of the
form $U_{d,\alpha_d}^{n_{d,\alpha_d}}-\theta_{l}U^{m^{(d,\alpha_d)}}$ are
irreducible in $\mathrm{gr}_{\nu}k_{(d-1)}[X_d]$. Thus the
decomposition above is the decomposition of $f$ into prime factors
in $\mathrm{gr}_{\nu}k_{(d-1)}[X_d]$. \par\noindent $(iii).$ It is
a result of $(iii)$ of last theorem and Lemma \ref{old-results},
$(iv)$.$\hfill\Box$

\begin{remark}{\rm Consider a valuation $\nu$ as above. The strong condition of Corollary \ref{irreducible elements}, $(ii)$,
is satisfied iff for any $i=0,\ldots,d-1$ either we have
$U_{i,\alpha_i}=0$ or $U_{i,\alpha_i}\neq 0$ and
$n_{i,\alpha_i}=\infty$. }
\end{remark}

\begin{theorem}\label{graded ring representation}{\bf (Homogeneous decomposition)} Let $\nu$ be a valuation associated to an {\rm SKP}. Consider the ring $R=k_{((\alpha,d))}$ and the restriction of $\nu$ to it. Every element $f\in R$ has a unique decomposition of the form $$f=p(T_{i_1},\ldots,T_{i_{d_1}})U^J,\ \mathrm{in}\ \mathrm{gr}_{\nu}R_{\nu},$$ where $d_1\leq d+1$ and  $A=\{i_1,\ldots,i_{d_1}\},$ for any $i\in A,$ $n_{i,\alpha_i}\neq \infty$ and $T_i=U_{i,\alpha_i}^{n_{i,\alpha_i}}U^{-m^{(i,\alpha_i)}}.$  And $0\leq J_{i',j}<n_{i',j},$ for $1\leq j\leq \alpha_{i'}.$ And $p(V_1,\ldots,V_{d_1})\in k[V_1,\ldots,V_{d_1}].$ \end{theorem} \noindent {\bf Proof.} This is a simple induction on $d$, using Theorem \ref{lemma representation graded}. For example if $U_{d,\alpha_d}\neq 0$, $n_{d,\alpha_d}\neq 0$ then $f=p(T_d)U_d^{J_d}$, in $\mathrm{gr}_{\nu}k_{(d-1)}[X_d]$, where the coefficients of $p(T_d)=\sum_lp_jT_d^l$ have the same $\nu$-value. By induction hypothesis, we have $p_l=q_l(T_{i_1},\ldots,T_{i_{d^l}})U_{d-1}^{J^l}$, in $\mathrm{gr}_{\nu}R_{\nu}$, where $d^l\leq d$. Now, all the $p_l$'s have the same $\nu$-value thus the vectors $J^l$ are the same for any $l$ (similar argument like proof of $(ii)$ of Theorem \ref{lemma representation graded}), we denote this vector by $J_{d-1}$. Hence, we have  $f=(\sum_lp_j)U_{d-1}^{J_{d-1}}U_{d}^{J_{d}}$ and we are done.    $\hfill\Box$

\begin{theorem}\label{omega type SKP}
Fix an {\rm SKP} $[U_{i,j},\beta_{i,j}]_{i=0..d,j=1..\alpha_i},$  such that $\alpha_d\geq \omega$. Suppose there
exists an infinite sequence of ordinals
$s_1<\cdots<s_{\omega}=\alpha_d$ such that $n_{d,s_j}>1,$ for any
$j<\omega$. Consider the acceptable vectors $\alpha^{(s_j)}$ (see
Definition \ref{delta-Euclid}). For any $f\in k_{(d-1)} [[X_d]]$
there exists  $j_{*}\in\mathbb{N}$ such that for any $j\geq j_{*}$
we have
$$\nu_{\alpha^{(s_j)}}(f)=\nu_{\alpha^{(s_{j_{*}})}}(f).$$
 Thus the limit $\lim_{j\to \omega}\nu_{\alpha^{(s_j)}}$ is
 well-defined and is equal to $\nu_{\alpha^{(s_{\omega})}}=\nu$.
\end{theorem}
\noindent {\bf Proof.} Multiplying $f$ by a suitable factor $u\in k^{(d-1)}$ we can assume $f\in k^{(d)}$. By assumptions, we have $U_{d,\alpha_d}=0$. Thus by Corollary \ref{one initial}, we have $\mathrm{in}_{\nu_{\alpha}}(f)=c_JU_d^{J}$, $c_J\in k^{(d-1)}$. Suppose  $j_*$ is the maximum index such that $J_{s_{j_*}}\neq 0$. Then by the algorithm of getting $adic$ expansion, this $j_{*}$ satisfies the conclusion of the Theorem. $\hfill\Box$

  }

\section{SKP-Valuations and numerical invariants }
{ One of the   ways to classify valuations is through their
numerical invariants. In this section we show how the arithmetic
of the SKP's  of an SKP-valuation  determines the
numerical invariants of the associated valuation on the field
$\ke^{(d)}$.

 We define the notion of pseudo-SKP. It allows us to avoid
ordinal numbers greater than $\omega$ for $\alpha_i$.

\begin{definition}\label{ pseudo-SKP}
{\rm For a {\rm SKP} $[U_{i,j},\beta_{i,j}]_{i=0..d,j=1..\alpha_i}$ a pseudo-SKP is a subset of $U$'s and $\beta$'s
which comes from dropping an arbitrary number of $U_{i,j}$'s (and
associated $\beta_{i,j}$'s) for $j<\alpha_i$ such that $n_{i,j}=1$. To any SKP is associated a minimal pseudo-SKP
which is obtained by dropping all $U_{i,j}$ for $j<\alpha_i$ such that $n_{i,j}=1$. This
minimal associated pseudo-SKP is unique. We denote this minimal
pseudo-SKP by $[U_{i,j},\beta_{i,j}]_{i=0..d,j=1..\alpha '_i}$, where $\alpha '_i\leq \omega$ (using the same notation as {\rm
SKP}'s). }
\end{definition}

\begin{proposition}

Fix an {\rm SKP} $[U_{i,j},\beta_{i,j}]_{i=0..d,j=1..\alpha_i}$ and let $\nu$ be the associated $k-valuation$.  Let
$[U_{i,j},\beta_{i,j}]_{i=0..d,j=1..\alpha '_i}$ be its
minimal pseudo-SKP. The valuation $\nu$ can be defined using the
data $[U_{i,j},\beta_{i,j}]_{i=0..d,j=1..\alpha '_i}$.
\end{proposition}

\noindent {\bf Proof.} To define the
valuation $\nu$ it is sufficient to know the $adic$ expansion of
elements. Moreover, in the $adic$ expansion of an element the
$U_{i,j}$'s with $n_{i,j}=1$ cannot appear. Thus the $adic$
expansion of every element is defined using only the  minimal
pseudo-SKP associated to $\nu$.$\hfill\Box$

The following lemma computes the rank and rational rank and value semigroup of an SKP valuation in terms of the
 arithmetic of the SKP.
\begin{lemma} \label{inductive-classification}
Consider a  $k-$valuation centered on the ring
$\ke^{(d)}$ such that
$\nu=\mathrm{val}[U_{i,j},\beta_{i,j}]_{i=0..d,j=1..\alpha_i}$. Let
$\overline{\nu}=\nu\mid_{k_{(d-1)}}.$ By Remark \ref{after def}.$(v)$ the data $[U_{i,j},\beta_{i,j}]_{i=0..d-1,j=1..\alpha_i}$ is an {\rm SKP}.
\begin{itemize}
\item[$(i)$]  We have $\overline{\nu}=\mathrm{val}[U_{i,j},\beta_{i,j}]_{i=0..d-1,j=1..\alpha_i}$.
\item[$(i)$] We have $\mathrm{rk}(\nu)-\mathrm{rk}(\overline{\nu})\in\{0,1\}.$ More precisely  $\mathrm{rk}(\nu)=\mathrm{rk}(\overline{\nu})+1$ iff
$\beta_{d,\alpha_d}\notin \Delta$ $(\Delta$ is the smallest isolated subgroup of $\Phi$ such that $\Phi^{*}_{d-1,\alpha_{d-1}}\subset \Delta),$ and $\mathrm{rk}(\nu)=\mathrm{rk}(\overline{\nu})$ iff $\beta_{d,\alpha_d}\in \Delta.$
\item[$(ii)$] We have $\mathrm{r.rk}(\nu)-\mathrm{r.rk}(\overline{\nu})\in\{0,1\}.$ More precisely  $\mathrm{r.rk}(\nu)=\mathrm{r.rk}(\overline{\nu})+1$ iff
$\beta_{d,\alpha_d}\notin \Phi^{*}_{d-1,\alpha_{d-1}},$ and
$\mathrm{r.rk}(\nu)=\mathrm{r.rk}(\overline{\nu})$ iff
$\beta_{d,\alpha_d}\in \Phi^{*}_{d-1,\alpha_{d-1}}.$
\item[$(iii)$] The semigroup $\nu(k^{(d)}\setminus\{0\})$ is equal to $\Gamma_{d,\alpha_d}$.
\end{itemize}
\end{lemma}

\noindent
{\bf Proof. } We only prove $(i)$. It is the consequence of the fact that for any $f\in k^{(d-1)} $ the $adic$ expansions of  $f$ with respect to the  two SKP's $[U_{i,j},\beta_{i,j}]_{i=0..d,j=1..\alpha_i}$  and     $[U_{i,j},\beta_{i,j}]_{i=0..d-1,j=1..\alpha_i}$ are the same. $\hfill\Box$

\begin{theorem}\label{classification}
Consider a  $k-$valuation centered on the ring $k[[X_0,X_1,X_2]],$
$\nu,$ which is defined by an {\rm SKP,} i.e., let
$\nu=\mathrm{val}[U_{i,j},\beta_{i,j}]_{i=0,1,2,j=1..\alpha_i}.$
Moreover, we suppose $\beta_{0,1}\in\Delta_1.$ Then we can compute
the numerical invariants of this valuation using the arithmetic of its minimal pseudo-{\rm SKP.}
This is summarized in
 {\rm Table} \ref{tabel-numeric}.
\begin{table}%[H]
\begin{center}
\hspace*{-2.8cm}
\begin{tabular}{p{6in}}
\begin{tabular}{| c| c| c| c| c| c| }
\hline
  &\multicolumn{2}{c}{Arithmetic of minimal  pseudo-SKP of the valuation $\nu$}&rk&r.rk&tr.deg\\
\hline

%%%%%%%%%%%%%%%%%%%%%%%%%%%%%%%%%%%%%%%%%%%%%%%%%%%%%%%%%%%%%%%%%%%%%%%%%%%%%%%%%%%%%%%%%%
$(\mathrm{I})$&$\alpha '_1<\infty,\ \alpha '_2<\infty $&$\beta_{i,j}\in\mathbb{Q}\beta_{0,1}$&$1$&$1$&$2$\\
 \hline
%%%%%%%%%%%%%%%%%%%%%%%%%%%%%%%%%%%%%%%%%%%%%%%%%%%%%%%%%%%%%%%%%%%%%%%%%%%%%%%%%%%%%%%%%%
\begin{tabular}{c}
$(\mathrm{II})_1$\\
  $(\mathrm{II})_2$
\end{tabular}&
\begin{tabular}{c}
$\alpha '_1<\infty,\ \alpha '_2<\infty $\\ $\alpha '_1<\infty,\ \alpha '_2<\infty $
\end{tabular}&
\begin{tabular}{c}
$\beta_{i,j}\in \Delta_1,\
\beta_{1,\alpha '_1}\in\mathbb{Q}\beta_{0,1},$
$\beta_{2,\alpha '_2}\in\Delta_1\backslash \mathbb{Q}\beta_{0,1}$\\
 $\beta_{i,j}\in  \Delta_1,\
\beta_{1,\alpha '_1}\in\Delta_1\backslash\mathbb{Q}\beta_{0,1},$
$\beta_{2,\alpha '_2}\in \mathbb{Q}\beta_{0,1}$
\end{tabular}
&$1$&$2$&$1$\\
\hline

%%%%%%%%%%%%%%%%%%%%%%%%%%%%%%%%%%%%%%%%%%%%%%%%%%%%%%%%%%%%%%%%%%%%%%%%%%%%%%%%%%%%%%%%%%%
 \begin{tabular}{c}
 $(\mathrm{III})_1$\\
 $(\mathrm{III})_2$
 \end{tabular}&
 \begin{tabular}{c}
 $\alpha '_1=\infty,\ \alpha '_2<\infty$\\
 $\alpha '_1<\infty,\ \alpha '_2=\infty$
 \end{tabular}&
\begin{tabular}{c}
$\beta_{i,j}\in\mathbb{Q}\beta_{0,1}$\\
$\beta_{i,j}\in\mathbb{Q}\beta_{0,1}$
\end{tabular}&
$1$&$1$&$1$\\\hline
%%%%%%%%%%%%%%%%%%%%%%%%%%%%%%%%%%%%%%%%%%%%%%%%%%%%%%%%%%%%%%%%%%%%%%%%%%%%%%%%%%%%%%%%%%%
(IV)&$\alpha '_1<\infty,\alpha '_2<\infty$&$
\beta_{1,\alpha '_1}\in\Delta_1\backslash\mathbb{Q}\beta_{0,1},\
\beta_{2,\alpha '_2}\in\Delta_1\backslash(\beta_{0,1},\beta_{1,\alpha '_1})\otimes\mathbb{Q}$&$1$&$3$&$0$\\
\hline
%%%%%%%%%%%%%%%%%%%%%%%%%%%%%%%%%%%%%%%%%%%%%%%%%%%%%%%%%%%%%%%%%%%%%%%%%%%%%%%%%%%%%%%%%%%%
\begin{tabular}{c}
$(\mathrm{V})_1$\\ $(\mathrm{V})_2$
\end{tabular}&
\begin{tabular}{c}
$\alpha '_1=\infty, \alpha '_2<\infty$\\  $\alpha '_1<\infty,\alpha '_2=\infty$
\end{tabular}&
\begin{tabular}{c}
$
\beta_{2,\alpha '_2}\in\Delta_1\backslash\mathbb{Q}\beta_{0,1} $\\
$ \beta_{1,\alpha '_1}\in\Delta_1\backslash\mathbb{Q}\beta_{0,1}$
\end{tabular}&
$1$&$2$&$0$\\\hline
%%%%%%%%%%%%%%%%%%%%%%%%%%%%%%%%%%%%%%%%%%%%%%%%%%%%%%%%%%%%%%%%%%%%%%%%%%%%%%%%%%%%%%%%%%%%
\begin{tabular}{c}
$(\mathrm{VI})$
\end{tabular}&
\begin{tabular}{c}
$\alpha '_1<\infty, \alpha '_2<\infty$
\end{tabular}&

 \begin{tabular}{c}$\mathrm{max}\{\beta_{i,\alpha '_i}\}\in\Delta_2\backslash\Delta_1,\
 \beta_{1,\alpha '_1}\in(\beta_{0,1},\beta_{2,\alpha '_2})\otimes\mathbb{Q}$
\end{tabular}&
$2$&$2$&$1$\\\hline

%%%%%%%%%%%%%%%%%%%%%%%%%%%%%%%%%%%%%%%%%%%%%%%%%%%%%%%%%%%%%%%%%%%%%%%%%%%%%%%%%%%%%%%%%%%%

\begin{tabular}{c}
$(\mathrm{VII})_1$\\ $(\mathrm{VII})_2$
\end{tabular}&
\begin{tabular}{c}
$\alpha '_1<\infty, \alpha '_2<\infty$\\ $\alpha '_1<\infty, \alpha '_2<\infty$
\end{tabular}&
\begin{tabular}{c}$
\beta_{1,\alpha '_1}\in\Delta_2\backslash\Delta_1,\
 \beta_{2,\alpha '_2}\in\Phi\backslash\Delta_2$\\
$\beta_{2,\alpha '_2}\in\Delta_2\backslash\Delta_1,\
 \beta_{1,\alpha '_1}\in\Phi\backslash\Delta_2$
\end{tabular}&
$3$&$3$&$0$\\\hline
%%%%%%%%%%%%%%%%%%%%%%%%%%%%%%%%%%%%%%%%%%%%%%%%%%%%%%%%%%%%%%%%%%%%%%%%%%%%%%%%%%%%%%%%%%%%
\begin{tabular}{c}
$(\mathrm{VIII})_1$\\
$(\mathrm{VIII})_2$
\end{tabular}&
\begin{tabular}{c}
$\alpha '_1=\infty, \alpha '_2<\infty$\\  $\alpha '_1<\infty,\alpha '_2=\infty$
\end{tabular}&
\begin{tabular}{c}$\beta_{2,\alpha '_2} \in\Delta_2\backslash\Delta_1$\\
 $\beta_{1,\alpha '_1}\in \Delta_{2}\backslash\Delta_1$\\

\end{tabular}&
$2$&$2$&$0$\\\hline
%%%%%%%%%%%%%%%%%%%%%%%%%%%%%%%%%%%%%%%%%%%%%%%%%%%%%%%%%%%%%%%%%%%%%%%%%%%%%%%%%%%%%%%%%%%%
\begin{tabular}{c}
$(\mathrm{IX})$
\end{tabular}&
\begin{tabular}{c}
$\alpha '_1<\infty, \alpha '_2<\infty$
\end{tabular}&
\begin{tabular}{c}$\mathrm{max}\{\beta_{i,\alpha '_i}\}\in\Delta_2\backslash\Delta_1,\
 \beta_{1,\alpha '_1}\in\Delta_2\backslash(\beta_{0,1},\beta_{2,\alpha '_2})\otimes\mathbb{Q}
 $
\end{tabular}&
$2$&$3$&$0$\\\hline
%%%%%%%%%%%%%%%%%%%%%%%%%%%%%%%%%%%%%%%%%%%%%%%%%%%%%%%%%%%%%%%%%%%%%%%%%%%%%%%%%%%%%%%%%%%%
(X)&$\alpha '_1=\infty,\alpha '_2=\infty$&&$1$&$1$&$0$\\\hline
%%%%%%%%%%%%%%%%%%%%%%%%%%%%%%%%%%%%%%%%%%%%%%%%%%%%%%%%%%%%%%%%%%%%%%%%%%%%%%%%%%%%%%%%%%%%
\end{tabular}
\end{tabular}
\end{center}
\caption{Numerical invariants via arithmetic of SKP of the
valuation  \label{tabel-numeric}}
\end{table}
\end{theorem}
%WHAT FOLLOWS IS NOT A PROOF. IT IS MUCH TOO SKETCHY
\noindent {\bf Proof.} The computation of the rank and the rational-rank
is a simple task. The only nontrivial task is the computation of
the transcendence degree or the dimension of valuation. It is a direct
calculation using Theorem \ref{graded ring representation}. For
example in the case $(I)$, pick $f,g\in k_{(d-1)}$ with
$\nu(f)=\nu(g)$. Then by Theorem \ref{graded ring representation}
we have $\mathrm{in}(f)=p(T_1,T_2)U^{J}$ and $\mathrm{in}(g)=q(T_1,T_2)U^{J'}$.
Using the properties of $J$ and $J'$ in the theorem, we see that
$J=J'$. Thus $f/g=p(T_1,T_2)/q(T_1,T_2)$. This shows
$k_{\nu}=R_{\nu}/\frak{m_{\nu}}=k(T_1,T_2)$.
We show  that $T_1$ and $T_2$ are algebraically
independent in $k_{\nu}.$ If $T_2$ is algebraic over $k(T_1)$, then there is a polynomial $0\neq p(T)\in k(T_1)[T]$ such that $p=p(T_2)=\sum_i c_iT_2^i=0$ in $k_{\nu}$. Regarding $T_1$ and $T_2$ as elements of $R_{\nu}$, we have $p(T_2)=\sum_i c_iT_2^i\in \frak{m}_{\nu}$. Note that $T_1=\frac{U_{1,\alpha_1}^{n_{1,\alpha_1}}}{U^{m^.}}$ and $T_2=\frac{U_{2,\alpha_2}^{n_{2,\alpha_2}}}{U^m}$.  Multiplying $p$ with a suitable power of  $U^{m^.+m}$, say n, we can assume that $U^{n(m^.+m)}p\in k_{((\alpha,2))}$.  The condition $p\in \frak{m}_{\nu}$ implies that the cancelation should occur between initail monomials of monomials of $U^{n(m^.+m)}p$ in the course of getting the $adic$ expansion. We show that this is impossible.\par Write $p=\sum_{i,j}r_{i,j}T_1^iT_2^j,$ $r_{i,j}\in k$. Then $U^{n(m^.+m)}p=\sum_{i,j}r_{i,j}U^{n_{[i,j]}}U_{1,\alpha_1}^iU_{2,\alpha_2}^j$. By Lemma \ref{tec}.$(ii)$ no cancelation can occur between initial monomials of  monomials of $U^{n(m^.+m)}p$  with different $j$'s (notice that index $(2,\alpha_2)$ does not occur in $U^{m_{[i,j]}}$). It remains to show that no cancelation can occur for a sum of the form $q_j=\sum_{i} r_{i,j}U^{n_{[i,j]}}U_{1,\alpha_1}^iU_{2,\alpha_2}^j$.  Notice that the power of $U_{2,}$'s,  are the same for different monomials of $q$ and the power of $U_{1,\alpha_1}$ are different for any two monomial of $q$. Now, the proof of  Lemma \ref{tec}.$(ii)$  shows that in the course of getting the $adic$ expansion of the monomials of $q$ the power of  $U_{1,\alpha_1}$ in the initial monomials  remain diffrent, for any two monomial of $q$.  Thus  no cancelation can occur between the initial monomials of $q$.  $\hfill\Box$
}

\section{Realization of a certain class of semi-groups as value semi-groups of polynomial rings }
{
  In this section we give a result on the  realization of a semi-group as the  semi-group of values which takes a valuation  on a polynomial ring.

% we define  $\mathrm{g.o.l}(\Gamma)$=n, and if $\alpha\geq \omega^2$ we define $\mathrm{g.o.l}(\Gamma)=\infty$. Thus, when $\mathrm{g.o.l}(\Gamma)<\infty$ the minimal system of generators of $\Gamma$ is of ordinal type $\omega\mathrm{g.o.l}(\Gamma)$.
\noindent
\begin{theorem}\label{semigroup generation}
 Let $\Gamma$ be a semigroup of an ordered abelian group   $(\Psi,<)$, given by a minimal system of generators $\{\gamma_j\}_{j\leq \alpha}\subseteq\Psi^+$, where  $\alpha=\omega n+j$, for $n,j\in \mathbb{N}$. Suppose  $\Gamma$ is of positive type (Definition \ref{positively generated}), and $\gamma_{j+1}> n_j\gamma_j$ when $n_j\neq \infty$. Set $G=(\Gamma)$ and $d=\mathrm{r.rk}(G)$.  Then there exists a zero-dimensional valuation $\nu$ of the field  $k(X_1,\ldots,X_d)$, centered on the polynomial ring $R=k[X_1,\ldots,X_d]$, such that its value-semigroup  is equal to $\Gamma$.
\end{theorem}
\noindent
{\bf Proof. } %We give special names for those indices of the $\gamma$'s that are not rationally dependent to the previous ones $\gamma_{s^t_{t'}}$. Introducing a new variable for every $\gamma_{s^t_{t'}+1}$, we construct  a set of key-polynomials of this new variable with values equal to $\gamma$, up to $\gamma_{s^t_{t'+1}}$. Then, we need to define a new variable. However, the situation differs in the case of limit ordinal:  If $s^t_{t'+1}$  is a limit ordinal then the limit key-polynomial which is available is zero and cannot take $\gamma_{s^t_{t'+1}}$ as its value. Thus, we are forced to define a new variable for taking value $\gamma_{s^t_{t'+1}}$. The precise definition is as follows.
Consider the semigroup $\Gamma$ with the minimal systems of generators  $\{\gamma_{j'}\}_{j'\leq \alpha}$, and suppose  $\alpha=\omega n+j^*$, $j^*\in\mathbb{N}$. We give new names ($s^t_{t'}$) to the indices of those $\gamma$'s which  are rationally independent from the previous ones (by  Proposition \ref{generalized effective component} this includes all the indices which are limit ordinals, i.e.,  for $t=1..n$ we have $n_{\omega t}=\infty$, see Lemma \ref{halgiri} for definition of $n$):  For $t=1..n+1$, let $f_t\in \mathbb{N}$ be the number of $j'$ such that $\omega(t-1)\leq j'<\omega t$ and $n_{j'}=\infty$, then  set $s^t_{t'}:=j'$ when $j'$ is the  $t'$th  such  $j'$,  for $t'=1..f_t$ . % We define $f_0:=0$, and  $s^{t}_{f_{t}+1}:=s^{t+1}_1$, for $t=1..n$.
Then we have  $$\{\gamma_{j'}\}_{j'\leq \alpha}=\{\gamma_{s^t_{t'}+j}\}_{t=1..n+1,t'=1..f_{t},j=0..j_{t,t'}},$$
 where $j_{t,t'}$ is the number of indices $j'$ such that $\gamma_{s^t_{t'}}\leq \gamma_{j'}<\gamma_{s^t_{t'+1}}$.  %$J_{t,t'}=\{1,\ldots,s^t_{t'+1}-s^t_{t'}\}$ for $t=1..n+1$ and $t'<f_t$, $J_{t,f_t}=\{j:\ 1\leq j<\omega\}$ for $t=1..n$,    $J_{n+1,f_{n+1}}=\{j:\ 1\leq j\leq j^*\}$ where $\alpha=\omega n+j^*$, $j^*\in\mathbb{N}$ (for the simplicity of notation we assume $f_{n+1}<j^*$), and finally $J_{t,f_{t}+1}=\{0\}$.
 %The finite set $\{\gamma_{s^t_{t'}}\}_{t'\in T'_t}$ are all those $\gamma$, $\gamma_{\omega(t-1)}\leq\gamma<\gamma_{\omega(t)}$, with $n_{s^t_{t'}}=\infty$ (see Lemma \ref{halgiri}, for definition of $n$).

 \par Then, by Proposition \ref{generalized effective component}, we have $\mathrm{r.rk}(G)=f_1+\cdots+f_{n+1}$.

\par   We define new indices $i_{t,t'}$ which will be the indices of the variables of the polynomial ring:  For $t=1..n+1$ and  $t'=1..f_t$ set $i_{t,t'}:=f_0+\cdots+f_{t-1}+t'$, where by convention $f_0=0$. %and $i_{t,t'}:=f_0+\cdots+f_t+t+1$, for $t'=f_{t}+1$,
  The total number $d$ of $i_{t,t'}$'s which has been defined is equal to:

\begin{displaymath}
d= %\left \{ \begin{array}{lll}
            f_1+\cdots+f_{n+1}%&
            =\mathrm{r.rk}(G).%+\mathrm{g.o.l}(\Gamma)  & :j^*=f_{n+1}\\
           % f_1+\cdots+f_{n+1}+1&=\mathrm{r.rk}(G)+\mathrm{g.o.l}(\Gamma)+1 & :j^*>f_{n+1}\\ \end{array} \right.
\end{displaymath}

  It is straightforward to check that the sequence $\{\beta_{i_{t,t'},j}:=\gamma_{s^t_{t'}+j-1}\}_{i_{t,t'}=1..d,j=1..j_{t,t'}}$ is a sequence of values (note the index $i$ starts from $1$). The key-polynomials of the  SKP associated to this sequence of values are elements of the ring $R=k[X_1,\ldots,X_{d}]$. The valuation $\nu$ associated to this SKP has value semi-group $\Gamma$. \par Notice that we have $\mathrm{r.rk}(\nu)=\mathrm{dim}R=d$. Hence, we are in the case of equality of Abhyankar's inequality $\mathrm{r.rk}(\nu)+\mathrm{tr.deg}(\nu)\leq \mathrm{dim}R=d$. Thus, the valuation $\nu$ is zero-dimensional. Moreover, one can not realize $\Gamma$ as a value semi-group of a polynomial ring with $<d$ variables.   $\hfill\Box$

 \begin{remark} {\rm The following remarks are in order:
 \begin{itemize}
 \item The positivity condition is quite restrictive in general. However, in the case we restrict to the value semi-groups of polynomial rings of two variables, all the value semi-groups are of positive type (See Proposition 4.2 of \cite{CutTe}). Moreover, in this case, if the ordinal type of the group is $\omega^2$ we are in the equality case of Abhyankar's inequality and the semigroup has to be of positive type
 %\item $\mathrm{tr.deg}(\nu)=0$. It seems that the bound $d$ obtained  for the number of the variables of the polynomial ring is the best possible. More precisely: we conjecture that given any semi-group $\Gamma$,  it cannot be realized as a value semigroup of a polynomial ring with $<d$ variables. %If this be true, as a result, we see that the minimal system of generators of the value semi-groups of polynomial rings has ordinal type $<\omega h$, for some $h\in\mathbb{N}$.
 \item The semigroup  $\Gamma$ is well ordered by \cite{Nu}, it is of ordinal type $\leq \omega^{\mathrm{rk}(G)}$ by (\cite{ZS} Vol.II, Appendix 3, Proposition 2).
 \end{itemize}
 }
 \end{remark}

 \begin{proposition}\label{generalized effective component}
 With the notation of Theorem \ref{semigroup generation} and Lemma \ref{halgiri}, for any limit ordinal $\omega(i+1)\leq\alpha$ we have   $\mathrm{rk}(G_{\omega(i+1)})=\mathrm{rk}(G_{\omega(i+1)^{-}})+1$. In particular, $n_{\omega(i+1)}=\infty$.
 \end{proposition}
 \noindent
 {\bf Proof. } We extend the notion of effective component to this situation. Consider an order  embedding $(\Phi,<)\subseteq(\mathbb{R}^n,<_{lex})$ such that $\Gamma\subseteq\mathbb{R}^n_{\geq_{les} 0}$. By definition, the effective component for the limit ordinal $\omega i$ is the first index $t\leq n$ such that $\#\{(\gamma_j)_t\}_{\omega i\leq j <\omega(i+1)}=\infty$. Like in the case of effective components, one can prove $t$ is well-defined. Note that $(\gamma_j)_{t'}=0$, for $t'<t$ and $j<\omega (i+1)$. Moreover, one can show that the content of Proposition \ref{canonical simplicity}.$(i)$ and $(ii)$ hold in this case. Suppose the effective component for $\omega i$ is $t$. Then an argument similar to the proof of Proposition \ref{canonical simplicity}.$(iii)$, shows that $(\gamma_j)_t\to+\infty\ (j\to\omega (i+1))$. But $\gamma_{\omega(i+1)}>_{lex}\gamma_{\omega i+j}$, for $j\in\mathbb{N}$. This is possible only if  $(\gamma_{\omega(i+1)})_{t'}>0$, for some $t'<t$. $\hfill\Box$
}

\bibliography{valbib}
\bibliographystyle{plain}
\end{document}